\documentclass[11pt]{article}
\usepackage{hyperref}
\usepackage{epsfig}
\usepackage{graphicx}
\usepackage{amsfonts} 
\usepackage{amssymb}  
\usepackage{mathrsfs} 
\usepackage{amsmath}\numberwithin{equation}{section}  
\usepackage{verbatim} 
\usepackage{color}
\usepackage{booktabs}
\usepackage{array}
\usepackage{multirow}
\usepackage{caption}
\usepackage{float} 
\usepackage{subfigure}
\usepackage{algpseudocode}
\usepackage{algorithmicx,algorithm}
\usepackage[title]{appendix}

\def\range{\mathrm{\bf range}}
\def\dom{\mathrm{\bf dom}}
\def\argmin{\mathrm{argmin}}
\def\prox{\mathrm{\bf prox}}
\def\bfone{\mathbf{1}}
\def\aff{\mathrm{\bf aff}}

\def\ri{\mathrm{\bf ri}}
\def\spa{\mathrm{\bf span}}

\begin{document}
\newtheorem{thm}{Theorem}[section]
\newtheorem{de}{Definition}[section]
\newtheorem{lem}[thm]{Lemma}
\newtheorem{coro}[thm]{Corollary}
\newtheorem{prop}[thm]{Proposition}
\newtheorem{alg}{Algorithm}[section]
\newtheorem{ex}{Example}[section]

\newcommand{\ifs}{\mbox{\rm if~}}
\newcommand{\otherwise}{\mbox{\rm otherwise~}}
\newcommand{\prf}{\par{\bf Proof. }}
\newcommand{\solve}{\par{\bf Solution. }}

\def\bc{\mathbb{C}}
\def\bk{\mathbb{K}}
\def\bn{\mathbb{N}}
\def\bq{\mathbb{Q}}
\def\br{\mathbb{R}}
\def\bt{\mathbb{T}}
\def\bx{\mathbb{X}}
\def\by{\mathbb{Y}}
\def\bz{\mathbb{Z}}
\def\range{\mathrm{\bf range}}
\def\dom{\mathrm{\bf dom}}
\def\argmin{\mathrm{argmin}}
\def\prox{\mathrm{\bf prox}}
\def\bfone{\mathbf{1}}
\newcommand{\md}{\mathscr{D}}
\newcommand{\me}{\mathscr{E}}
\newcommand{\ms}{\mathscr{S}}
\newcommand{\mx}{\mathscr{X}}
\newcommand{\my}{\mathscr{Y}}
\newcommand{\mz}{\mathscr{Z}}
\newcommand{\bigvarphi}{\mbox{\Large $\varphi$}}
\newcommand{\supp}{{\rm supp}}
\newcommand{\cov}{\mathrm{Cov}}
\newcommand{\dist}{{\rm dist}}
\newcommand{\diam}{{\rm diam}}
\newcommand{\diag}{{\rm diag}}
\newcommand{\imf}{{\rm imf}}
\newcommand{\rank}{{\rm rank}}
\newcommand{\sgn}{{\rm sgn}}
\newcommand{\spann}{{\rm span}}
\newcommand{\re}{{\rm Re}}
\newcommand{\im}{{\rm Im}}
\newcommand{\aew}{{\rm ~a. e.~}}
\newcommand{\as}{{\rm ~a. s.~}}
\newcommand{\lip}{{\rm Lip}}
\newcommand{\bbox}{\hfill\rule{0.4em}{0.6em}}
\def\argmin{\mathop{\mathgroup\symoperators argmin}}
\def\argmax{\mathop{\mathgroup\symoperators argmax}}
\def\essup{\mathop{\mathgroup\symoperators essup}}
\newcommand{\embed}{\mathop{\subset\hspace{-1.5ex}_\succ}}
\def\st{\textrm{subject to~}}
\def\deg{\textrm{deg}}

\newcommand{\overcirc}[1]{\overset{\circ}{#1}}
\newcommand{\note}[1]{\textcolor{red}{\scriptsize(#1)}}

\renewcommand{\liminf}{\mathop{\mathgroup\symoperators\underline{lim}}}
\renewcommand{\limsup}{\mathop{\mathgroup\symoperators\overline{lim}}}

\newcommand{\bm}[1]{{\boldmath{\textbf{$#1$}}}}
\def\bfa{\bm{a}}
\def\bfb{\bm{b}}
\def\bfc{\bm{c}}
\def\bfe{\bm{e}}
\def\bff{\bm{f}}
\def\bfh{\bm{h}}
\def\bfn{\bm{n}}
\def\bfp{\bm{p}}
\def\bfq{\bm{q}}
\def\bfw{\bm{w}}
\def\bfu{\bm{u}}
\def\bfv{\bm{v}}
\def\bfx{\bm{x}}
\def\bfy{\bm{y}}
\def\bfz{\bm{z}}
\def\bfzero{\bm{0}}
\def\bftheta{\bm{\theta}}

\title{A proximal algorithm incorporating difference of convex functions optimization
for solving a class of single-ratio fractional programming
}
\author{Anna Qi\textsuperscript{a},
Jianfeng Huang\textsuperscript{b},
Lihua Yang\textsuperscript{c},
Chao Huang\textsuperscript{d}\footnote{Corresponding author. Email: hchao@szu.edu.cn}\\
\textsuperscript{a}{\small School of Mathematics and Systems Science, Guangdong Polytechnic Normal University, China} \\
\textsuperscript{b}{\small School of Financial Mathematics and Statistics, Guangdong University of Finance, China} \\
\textsuperscript{c}{\small School of Mathematics, Sun Yat-sen University, China}\\
\textsuperscript{d}{\small School of Mathematical Sciences, Shenzhen University, China}
}
\date{}
\maketitle
\begin{abstract}
\par
In this paper, we consider a class of single-ratio
fractional minimization problems, where both the numerator and denominator
of the objective are convex functions satisfying positive homogeneity.
Many nonsmooth optimization problems on the sphere that are commonly
encountered in application scenarios across different scientific fields can be
converted into this equivalent fractional programming.
We derive local and global optimality conditions of the problem
and subsequently propose a proximal-subgradient-difference of convex functions
algorithm (PS-DCA) to compute its critical points.
When the DCA step is removed, PS-DCA reduces to the proximal-subgradient
algorithm (PSA).
Under mild assumptions regarding the algorithm parameters, it is shown that any
accumulation point of the sequence produced by PS-DCA or PSA is a critical point
of the problem.
Moreover, for a typical class of generalized graph Fourier
mode problems, we establish global convergence of the entire sequence generated by
PS-DCA or PSA.
Numerical experiments conducted on computing the generalized graph Fourier
modes demonstrate that,
compared to proximal gradient-type algorithms, PS-DCA integrates difference of
convex functions (d.c.) optimization, rendering it less sensitive to initial
points and preventing the sequence it generates from being trapped in low-quality
local minimizers.

{\bf key words:} Single-ratio fractional programming, local and global optimality
conditions, proximal algorithm, d.c. optimization,
generalized graph Fourier mode problem
\end{abstract}

\newpage
\section{Introduction}

\setcounter{equation}{0}
A fundamental problem that emerges from applications in the fields of signal
processing and machine learning can be formulated as minimizing a convex
function with positive homogeneity on a unit sphere defined by a specific norm.
For instance, the graph Fourier basis problem\cite{2018irregular,dirGFT2017} and
the continuous relaxation of the Cheeger cut problem\cite{Rayleigh2019,IPM2010}.
Moreover, in tasks like dictionary learning and robust subspace recovery,
one often encounters the need to minimize a positively homogeneous convex
function over $\ell_2$ or $\ell_\infty$ spheres, e.g., see
\cite{Bai2019,dictlear2017,RSRdual2018,RSRalg2018} and the references therein.
We classify this type of optimization problem as the following model:
\begin{equation}
\label{eq:P1}
\min_{x\in X}T(x),~~~\text{s.t.}~~B(x)=1,
\end{equation}
where $X$ is a given linear space, $T(x)$ and $B(x)$ are convex functions taking
finite values on $\br^n$, and satisfy positive homogeneity and absolute
homogeneity across the entire domain, respectively:
\begin{equation}\label{eq:homogenuity}
\begin{cases}
T(\alpha x)=\alpha T(x)\ge 0,~~\forall\alpha\geq0,x\in \br^n,\\
B(\alpha x)=|\alpha|B(x)\ge 0,~~\forall\alpha\in\br,x\in \br^n.
\end{cases}
\end{equation}
Furthermore, it is also satisfied that $T(x)=0,~\forall x\in X_0$, where $X_0:=\{x\in X\mid B(x)=0\}$.

%
To gain a deeper understanding of how problem \eqref{eq:P1} is derived from
practical application problems, we explore one of its typical examples:
the graph Fourier basis problem,
which plays a central role in the emerging field
of signal processing on graphs\cite{graphon2022,harmonic2023,emerging2013,shuman2016}.
Various types of graph Fourier basis
can be defined by iteratively solving a sequence of graph signal variation minimization
problems\cite{1997Spectral,2018irregular,hammond2011,spreGFT2018,l1GFT2021},
such as the classical Laplacian basis, which can be formulated
as the solution to a sequence of $l_2$ norm minimization problems,
and the set of generalized graph Fourier modes, which is defined through
a sequence of graph variation minimization problems,
and the $k$-th graph Fourier mode is subject to $U_{k-1}^\top Qx=0$ and
$x^\top Qx=1$,
where $U_{k-1}$ is a matrix composed of the previous $k-1$ graph Fourier modes,
$Q$ is a preselected symmetric positive definite matrix. Such relaxation of the
inner product on the space of graph signals can take better account of the
irregularity of the graph structure\cite{2018irregular}.
In particular, when $Q$ is a diagonal positive definite matrix,
the resulting graph Fourier transform will yield desirable properties in physical sensing applications
\cite[\S\uppercase\expandafter{\romannumeral2}]{2018irregular}.
Slightly earlier than the introduction of the generalized graph Fourier modes,
S. Sardellitti et al. proposed to define the Fourier basis for directed graphs
based on graph directed total variation minimization to overcome the obstacle
that the classical Laplacian basis cannot be defined on directed
graphs\cite{dirGFT2017}. It is shown that the graph directed variation is
the exact continuous convex relaxation of the graph cut size
\cite[\S\uppercase\expandafter{\romannumeral2}]{dirGFT2017},
and consequently, minimizing the graph directed variation amounts to minimizing
the cut size function.
Although minimizing the graph directed total variation is generally not
equivalent to solving a sequence of graph directed variation minimization
problems iteratively, it can be proved that the solution of the latter is
also a solution of the former\cite[Property 1]{2018irregular}, thereby constituting the directed graph Fourier
basis.
Recently, a definition of graph Fourier basis based on
iteratively solving a sequence of $l_p$ norm variation minimization problems
was proposed, where the $k$-th basis vector satisfies the
the orthonormal condition: $U_{k-1}^\top x=0$ and $x^\top x=1$\cite{l1GFT2021}.
It is easy to verify that when $p\geq1$, the $l_p$ Fourier basis is the solution
of problem \eqref{eq:P1} if $T(x)=\Big(\sum_{1\leq i<j\leq n}w_{ij}|x_i-x_j|^p\Big)^{\frac1p}$,
$B(x)$ is the $l_2$ norm and the linear space $X=\{x\in\br^n|U_{k-1}^\top x=0\}$.
When we take $B(x)=\|Q^{\frac12}x\|$, where $Q$ is a symmetric positive definite
matrix, and let
$T(x)$ be any convex function satisfying positive homogeneity
(e.g. graph directed variation, $l_p$ norm variation and graph total variation etc),
the solution of problem \eqref{eq:P1} defines generalized graph Fourier modes.

Under the above assumptions, the problem \eqref{eq:P1} can be transformed into
solving a fractional programming
\begin{equation}\label{eq:P2}
\min_{x\in\br^m\setminus\{0\}}\frac{T(x)}{B(x)},
\end{equation}
where $T, B$ are convex functions taking finite values on $\br^m(m\leq n)$ and
satisfy condition \eqref{eq:homogenuity} on $\br^m$, furthermore,
$B(x)>0,~\forall x\in \br^m\setminus\{0\}$. For a rigorous proof, please refer
to the Appendix \ref{appendix}, which aims to establish the equivalence between
problem \eqref{eq:P1} and fractional minimization problem \eqref{eq:P2}.

In this paper, we
consider a fractional programming that is more general than \eqref{eq:P2},
as detailed below:
\begin{equation*}\label{eq:P}
P_1:~~  \min_{x\in\Omega}\frac{T(x)}{B(x)},
\end{equation*}
where $T(x)$ and $B(x)$ are convex functions taking finite values on $\br^n$, so they are both continuous on $\br^n$\cite[Th10.1]{1997Convex}, and both of them satisfy positive homogeneity,
\begin{equation}\label{eq:positive homogenuity}
\begin{cases}
T(\alpha x)=\alpha T(x)\ge 0,~~\forall\alpha\geq0,x\in\br^n,\\
B(\alpha x)=\alpha B(x)\ge 0,~~\forall\alpha\geq0,x\in\br^n.
\end{cases}
\end{equation}
Furthermore, we also assume that $\{x\in\br^n\mid B(x)=0\}\cap\{x\in\br^n\mid T(x)=0\}=\{0\}$ and the set $\Omega:=\{x\in\br^n\mid B(x)\neq0\}$ is nonempty.
Obviously, Problem $P_1$ includes \eqref{eq:P2} as a special case.
Furthermore, Problem $P_1$ inherently subsumes the well-known
generalized Rayleigh quotient problem\cite{Robust2018,Rayleigh2019}.

The research on fractional programming has a longstanding history,
the classic algorithm for solving a generalized single-ratio fractional
minimization problem is the parametric approach, e.g., the Dinkelbach type
algorithm and its modified versions, which solves a parametric suboptimal
problem in each loop to approximate the solution of the fractional
programming\cite{Dinkelbach1967,Parametric1983,modDink1975}. Crouzeix et al.
proved that if there is only one optimal
solution for the fractional programming and all parametric subproblems have
solutions, then the Dinkelbach type algorithm will linearly converge to the
optimal value of the fractional programming\cite{Crouzeix1985}. Although the parametric
approach possesses elegant theory, it is difficult to be used in practice when
applied to Problem $P_1$. In fact, the subproblem with respect to parameter $\lambda>0$
at each iteration takes the form of
\begin{equation*}\label{param approach}
\min_{x\in\Omega}\{T(x)-\lambda B(x)\},
\end{equation*}
which is essentially a d.c. optimization and generally nonconvex
(possibly nonsmooth), thus numerical computation is difficult.
Many efficient methods based on parametric approach usually assumed that the
parametric subproblem is convex to overcome this obstacle
\cite{generFP1996,dualmethod2017,inexapromx2008}.

Over the past few years, proximal-gradient type algorithms have been
extensively studied for addressing a wide variety of fractional
programming problems(see
\cite{proxgrad2017,singleloop2025,PGSABE2022,bergmanprox2024,PGSA2022} and
the references therein). In \cite{proxgrad2017}, the authors consider a single-ratio
fractional minimization
problem over a closed convex set, where the numerator is convex and the
denominator is required to be smooth, either concave or convex, the proposed
method can be applied to Problem $P_1$ when $B$ is smooth and the feasible
region is convex and closed, the resulting algorithm computes the new iterate by
\begin{equation}\label{prox-grad}
x^{k+1}=\argmin_{x\in\Omega}\Big\{T(x)+\frac1{2\eta_k}\|x-(x^k+
c_k\eta_k\nabla B(x^k))\|^2\Big\},
\end{equation}
where $\eta_k>0,~c_k=\frac{T(x^k)}{B(x^k)}$. A proximity-gradient-subgradient
algorithm (PGSA) was proposed sequentially in \cite{PGSA2022} for solving a
different class of fractional programming which takes the form of
\begin{equation}\label{structured FP}
\min_{x\in\Omega}\Big\{\frac{f(x)+h(x)}{g(x)}\Big\},
\end{equation}
where $f$ is proper and continuous within its domain (possibly nonsmooth and
nonconvex), $g$ is convex, $h$ is Lipschitz differentiable with a positive
Lipschitz constant $L$, and $\Omega:=\{x\in\br^n\mid g(x)\neq0\}$ is nonempty.
Starting from an initial point $x^0$, PGSA generates a sequence of iterates by
\begin{equation*}
x^{k+1}\in\argmin_{x\in\br^n}\Big\{f(x)+\frac1{2\lambda_k}
\|x-(x^k-\lambda_k\nabla h(x^k)+\lambda_kc_ky^k)\|^2\Big\},
\end{equation*}
where $\{\lambda_k\}$ is an input sequence satisfying $0<\lambda_k<\frac1L$,
$y^k\in\partial g(x^k)$ and $c_k=\frac{f(x^k)+h(x^k)}{g(x^k)}$.
Problem $P_1$ is included in \eqref{structured FP} by taking
$f(x)=T(x)$, $g(x)=B(x)$ and $h(x)=0$, the resulting iterative scheme is
\begin{equation}\label{pgsa}
x^{k+1}\in\argmin_{x\in\br^n}\Big\{T(x)+\frac1{2\lambda_k}
\|x-(x^k+\lambda_kc_ky^k)\|^2\Big\}.
\end{equation}
Subsequential convergence of PGSA is also guaranteed under the assumption
that $f+h$ and $g$ do not attain $0$ simultaneously, which is not satisfied by
Problem $P_1$, to ensure any limit point of $\{x^k\}$ resides on $\Omega$.
Additionally, it has been further shown that
the entire sequence generated by PGSA converges to a critical point of the
objective function under suitable assumption that the objective function
is level-bounded, which is also not fulfilled by the case of Problem $P_1$
since $T$ and $B$ are positive homogeneous, the level set associated with the
objective function is not bounded. Recently, PGSA combined with backtracked
extrapolation (PGSA\_BE) was proposed in \cite{PGSABE2022} to possibly accelerate
PGSA for solving problems in the form of \eqref{structured FP}. However, one of
the underlying assumptions of the model is that the numerator and denominator
of the objective function cannot be zero simultaneously. Furthermore, it was
also assumed that the objective function is level-bounded in order to
conduct convergence analysis of the proposed algorithm.
Therefore, the convergence analysis of the above algorithms can not be
applied to Problem $P_1$.

In this paper, we derive local and global optimality conditions for Problem $P_1$
based on d.c. optimization and then propose to perform a one-step
d.c. algorithm (DCA) after the proximal-subgradient procedure in each iteration,
the proposed algorithm is called proximal-subgradient-d.c. algorithm (PS-DCA).
If the DCA step is eliminated, PS-DCA reduces to proximal-subgradient
algorithm (PSA),
whose iteration scheme coincides with PGSA except that
PSA adds an extra step in each iteration to normalize the iterate by the
$l_2$ norm.
Under the assumption of bounded algorithm parameters, we demonstrate that
PS-DCA or PSA exhibits monotonicity and compactness, leading us to establish
subsequential convergence of the algorithm.
Additionally, we deduce a closedness property of PSA with respect to
local minimizers of the objective function.
As an application of the proposed methods, we study the global sequential
convergence of PS-DCA or PSA for the generalized graph
Fourier mode problem based on minimizing the graph directed variation. Finally, the
numerical results are presented to illustrate the advantages of the proposed
method.

The rest of the paper is organized as follows. In Section \ref{prelim},
we introduce the notations employed in this paper and provide necessary
preliminaries on d.c. optimization. In Section \ref{optimal cond}, we derive
local and global optimality conditions for Problem $P_1$. In Section \ref{alg-PS-DCA},
we present an equivalent condition for critical points of Problem $P_1$ and then
propose our PS-DCA.
In Section \ref{converg-of-PSDCA}, we present PSA and prove subsequential
convergence property of PS-DCA and PSA. Furthermore,
we also prove a closedness property of PSA concerning local minimizers
of the objective function.
In Section \ref{apply to GFMs prob}, we establish global convergence property of
PS-DCA and PSA for a typical class of generalized graph
Fourier mode problems and then we offer a subgradient selection strategy within
the DCA procedure for this type of problem.
Numerical results are reported in Section \ref{exp}.
Section \ref{conclud} is a final conclusion.

\section{Notations and preliminaries}\label{prelim}
In this section, we start with notations that will be used throughout the paper
and then we shall introduce some fundamental preliminaries on d.c. optimization.
The dot product in $\br^n$ is defined as
$$\langle x,y\rangle:=x_1y_1+x_2y_2+\cdots+x_ny_n.
$$
The Euclidean norm, induced by the dot product, corresponds to the $l_2$ norm
$\|x\|=\langle x,x\rangle^{\frac12}$. The $l_1$ norm on $\br^n$ is given by
$\|x\|_1=\sum_{i=1}^n|x_i|$. The symbol $\Gamma(\br^n)$ denotes the set of all proper lower semicontinuous convex functions on $\br^n$. The conjugate function $f^*$ of $f\in\Gamma(\br^n)$ is a function belonging to $\Gamma(\br^n)$, which is defined by
\begin{equation}\label{defini conjugate}
f^*(y):=\sup_{x\in\br^n}(x^\top y-f(x)).
\end{equation}
For any $C\subseteq\br^n$, $\mathrm{\bf int}C$ denotes the interior of $C$, i.e.,
$$\mathrm{\bf int}C:=\{x\mid\exists \epsilon>0~\text{s.t.}~\mathcal{B}(x,\epsilon)\subseteq C\},
$$
where $\mathcal{B}(x,\epsilon):=\{y\mid\|y-x\|\leq\epsilon\}$ is the
ball centered at $x$ with radius $\epsilon>0$. The relative interior of $C\subseteq\br^n$ is defined as
$$\ri C:=\{x\in \aff C\mid\exists \epsilon>0~\text{s.t.}~\mathcal{B}(x,\epsilon)\cap(\aff C)\subseteq C\},
$$
where $\aff C$ denotes the affine hull of $C$.
For a function $f:\br^n\rightarrow\br\cup\{+\infty\}$, we denote by
$\dom f=\{x\mid f(x)<+\infty\}$ the effective domain of $f$.
For $x\in \dom f$, the subdifferential of $f$ at $x$ is defined by
$$\partial f(x):=\{v\mid f(y)\ge f(x)+v^{\top}(y-x),~~\forall y\in\br^n\}.$$
The element of $\partial f(x)$ is called the subgradient of $f$ at $x$. It is easy to see that $\partial f(x)$ is a closed convex set.
$\dom\partial f=\{x\mid\partial f(x)\neq\emptyset\}$ denotes the effective domain of subdifferential mapping $\partial f$.
Let $\rho\geq0$, one says the function $f:\br^n\rightarrow\br\cup\{+\infty\}$
is $\rho$-convex if
$$f(\alpha x+(1-\alpha)y)\leq\alpha f(x)+(1-\alpha)f(y)-\frac{\alpha(1-\alpha)\rho}{2}\|x-y\|^2, \forall\alpha\in (0,1),\forall x,y\in\br^n.$$
It has been shown that $f$ is $\rho$-convex if and only if
$f-\frac{\rho}{2}\|\cdot\|^2$ is convex on $\br^n$\cite{Castellani1998}.
The modulus of strong convexity of $f$ on $\br^n$ is defined as
\begin{equation}\label{strongconvmu}
\rho(f):=\sup\{\rho\geq0\mid f-\frac\rho2\|\cdot\|^2~\text{is convex on}~\br^n\}.
\end{equation}
If $\rho(f)>0$, $f$ is called strongly convex on $\br^n$.

The remaining part of this section is devoted to preliminaries on d.c.
optimization.
We follow the convention that $+\infty-(+\infty)=+\infty$
\cite{HiriartUrruty1989,dualDCA1988}.
Generally, a d.c. program is that of the form
\begin{equation}\label{genDC}
\alpha=\inf_{x\in\br^n}f(x):=f_1(x)-f_2(x),
\end{equation}
where $f_1,f_2\in\Gamma(\br^n)$.
It's usually assumed in a d.c. programming that $\dom(f_1)\subseteq\dom(f_2)$, which implies $\alpha$ is finite. Such a function $f$ is called a d.c. function on $\br^n$, and $f_1$, $f_2$ are called its d.c. components. By using the definition of conjugate functions, we have
\begin{align*}
\alpha &=\inf_{x\in\br^n}\{f_1(x)-\sup_{y\in\br^n}\{x^\top y-f^*_2(y)\}\}\\
&=\inf_{x\in\br^n}\inf_{y\in\br^n}\{f_1(x)-x^\top y+f^*_2(y)\}\\
&=\inf_{y\in\br^n}\{f^*_2(y)-\sup_{x\in\br^n}\{x^\top y-f_1(x)\}\\
&=\inf_{y\in\br^n}\{f^*_2(y)-f^*_1(y)\}.
\end{align*}
The d.c. program
\begin{equation}\label{dual dc}
\alpha=\inf_{y\in\br^n}\{f^*_2(y)-f^*_1(y)\}
\end{equation}
is called the dual program of \eqref{genDC}. Let $\mathcal{P}$ and
$\mathcal{D}$ denote the solution sets of problems \eqref{genDC} and
\eqref{dual dc}, respectively. The relationship between primal and dual
solutions is $\cup\{\partial f_2(x)|x\in\mathcal{P}\}\subseteq\mathcal{D}$
and $\cup\{\partial f^*_1(y)|y\in\mathcal{D}\}\subseteq\mathcal{P}$
(see, for example, \cite{dualDCA1988}). Thus, when solving the primal d.c. program, its dual program is also solved simultaneously and vice versa.

Based on duality and local optimality conditions, T.Pham Dinh et al.
proposed a d.c. optimization algorithm (d.c. algorithm, DCA) in
\cite{thirtydevepDCA2018,Pham1998}
for solving a d.c. program, the DCA is described as follows:
Given $x^0\in\dom f_1$, for $k=0,1,\cdots$, let
$$y^k\in\argmin_{y\in\partial f_2(x^k)}(f_2^*(y)-f_1^*(y)),$$
$$x^{k+1}\in\argmin_{x\in\partial f_1^*(y^k)}(f_1(x)-f_2(x)).$$
Note that for any $f\in\Gamma(\br^n)$ and any $v\in\br^n$, the following
three conditions: $x\in\partial f(v)$, $v\in\partial f^*(x)$ and
$f^*(x)+f(v)=v^\top x$ are equivalent to each other\cite[Th23.5]{1997Convex},
so the iterative scheme of DCA is equivalent to the following two convex maximization problems
\begin{equation}\label{compDCAy}
y^k\in\argmin_{y\in\partial f_2(x^k)}(y^\top x^k-f_1^*(y)),
\end{equation}
and
\begin{equation}\label{compDCAx}
x^{k+1}\in\argmin_{x\in\partial f_1^*(y^k)}(x^\top y^k-f_2(x)).
\end{equation}
respectively. However, convex maximization problems, as a special case of d.c.
programming, still remain nonconvex and are difficult to solve in general. In
practice, one often performs simplified form of DCA, which takes a point within the feasible region instead of solving these subproblems. More precisely, simplified DCA generates two sequences $\{x^k\}$ and $\{y^k\}$ by setting
\begin{equation}\label{generDCA}
y^k\in\partial f_2(x^k), ~~x^{k+1}\in\partial f_1^*(y^k).
\end{equation}
The convergence analysis of the simplfied DCA
is conducted in \cite[\S3.4]{Pham1998}, we review the definition of
critical points of a d.c. function and the basic convergence theorem below.
\begin{de}
A point $x^*$ is said to be a critical point of $f_1-f_2$ if $\partial f_1(x^*)\cap\partial f_2(x^*)\neq\emptyset$.
\end{de}
\begin{thm}\label{simp dca converg thm}
Suppose that the sequence $\{x^k\}$ and $\{y^k\}$ are generated by the simplified DCA. Then we have
\begin{equation}\label{dc-converg-Th}
\begin{split}
(f_1-f_2)(x^{k+1})&\leq(f_2^*-f_1^*)(y^k)
-\max\big\{\frac{\rho_2}2\|dx^k\|^2,\frac{\rho_2^*}2\|dy^k\|^2\big\}\leq
(f_1-f_2)(x^k)\\
&-\max\big\{\frac{\rho_1+\rho_2}2\|dx^k\|^2,\frac{\rho_1^*}2\|dy^{k-1}\|^2
+\frac{\rho_2}2\|dx^k\|^2,\frac{\rho_1^*}2\|dy^{k-1}\|^2
+\frac{\rho_2^*}2\|dy^k\|^2\big\},
\end{split}
\end{equation}
where $0\leq\rho_i<\rho(f_i)$ $(resp.~0\leq\rho^*_i<\rho(f^*_i)),(i=1,2)$, and
$\rho_i=0$ $(resp. ~\rho^*_i=0)$ if $\rho(f_i)=0$ $(resp. ~\rho(f_i^*)=0)$.
Specifically, if the supremum $\rho(f_i)$ $(resp. ~\rho(f_i^*))$ is attainable,
then $\rho_i$ $(resp. ~\rho^*_i)$ may take this value.
$dx^k:=x^{k+1}-x^k$ and $dy^k:=y^{k+1}-y^k$.
The equality $(f_1-f_2)(x^{k+1})=(f_1-f_2)(x^k)$ holds if and only if
$x^k\in\partial {f_1}^*(y^k),y^k\in\partial f_2(x^{k+1})$ and
$(\rho_1+\rho_2)dx^k=\rho^*_1dy^{k-1}=\rho^*_2dy^{k}=0$. In this case
\begin{itemize}
  \item $(f_1-f_2)(x^{k+1})=(f_2^*-f_1^*)(y^k)$ and $x^k,x^{k+1}$ are
  critical points of $f_1-f_2$ satisfying $y^k\in(\partial{f_1}(x^k)
  \cap\partial{f_2}(x^k))$ and $y^k\in(\partial{f_1}(x^{k+1})\cap\partial{f_2}(x^{k+1}))$,
  \item $y^k$ is a critical point of $f_2^*-f_1^*$ satisfying $[x^k,x^{k+1}]\subseteq
  (\partial {f_1^*}(y^k)\cap\partial {f_2^*}(y^k))$,
  \item $x^{k+1}=x^k$ if $\rho(f_1)+\rho(f_2)>0$, $y^k=y^{k-1}$ if
  $\rho(f_1^*)>0$ and $y^k=y^{k+1}$ if $\rho(f_2^*)>0$.
\end{itemize}
\end{thm}
The authors further point out that for the case of complete DCA,
when $(f_1-f_2)(x^{k+1})=(f_1-f_2)(x^k)$,
$x^k$ must satisfy $\partial {f_2}(x^k)\subseteq\partial {f_1}(x^k)$.

Theorem \ref{simp dca converg thm} indicates that the value of the strong
convexity modulus in both primal and dual problems can influence the
efficiency of DCA. Therefore, when applying the simplified DCA, one usually lets
$\tilde{f_1}=f_1+\frac\rho2\|\cdot\|^2$,
$\tilde{f_2}=f_2+\frac\rho2\|\cdot\|^2$,
where $\rho>0$ represents the regularization parameter,
which is chosen appropriately according to the d.c. decomposition adopted
by the objective function.
Additionally, the d.c. components in primal problem are all strongly convex
functions, which guarantees that their conjugate functions
are differentiable throughout the interior of
their respective effective domains\cite[Th25.1,Th26.3]{1997Convex}, i.e.,
$\partial\tilde{f_i}^*(x)=\{\nabla\tilde{f_i}^*(x)\},\forall x\in\mathrm{\bf int}(\dom f_i),(i=1,2)$. Consequently, if $\dom\partial\tilde{f_1}^*\subseteq
\mathrm{\bf int}(\dom f_1)$, solving \eqref{compDCAx} is equivalent to selecting $x^{k+1}$ from $\partial\tilde{f_1}^*(y^k)$, which is the same as simplified DCA.

\section{Local and global optimality conditions for Problem $P_1$}
\label{optimal cond}
In this section, we deduce local and global optimality conditions
for Problem $P_1$ based on d.c. optimization.
Denote $E(x):=\frac{T(x)}{B(x)}$, as $E$ is invariant under positive scalings,
Problem $P_1$ is further equivalent to Problem $P_2$,
\begin{equation}\label{eq:P4}
P_2:~~\min_{x\in\mathcal{S}\cap\Omega}E(x),
\end{equation}
where $\mathcal{S}=\{x\mid\|x\|=1\}$. Thus the discussion of local or global
minimizers of $E$ can be conducted on the set $ \mathcal{S}\cap\Omega$. To
this end, define the neighborhood of $x^*$ on $\mathcal{S}\cap\Omega$ by
$\mathcal{B}_{\scriptscriptstyle\mathcal{S}\cap\Omega}(x^*,\epsilon):=
\{x\in\mathcal{S}\cap\Omega\mid\|x-x^*\|\leq\epsilon\}$. With this notation,
we introduce the following definition of local minimizers.
\begin{de}\label{local min}
$x^*\in\mathcal{S}\cap\Omega$ is a local minimizer of $E$ if there exists
$\epsilon>0$ such that for $\forall x\in\mathcal{B}_
{\scriptscriptstyle\mathcal{S}\cap\Omega}(x^*,\epsilon)$,
there holds $E(x)\geq E(x^*)$. Furthermore, if  $E(x)>E(x^*),\forall x\in
\mathcal{B}_{\scriptscriptstyle\mathcal{S}\cap\Omega}(x^*,\epsilon)$ and
$x\neq x^*$, we say
$x^*$ is a strict local minimizer of $E$.
\end{de}
Let
\begin{equation}\label{define-gh}
g(x):=T(x), ~~~h(x):=E(x^*)B(x).
\end{equation}
Below, we present a proposition that elucidates
the relationship between the local minimizers of $E$ and those of the
d.c. function $g-h$.
\begin{prop}
\label{loc trans prop}
Assume $x^*\in\mathcal{S}\cap\Omega$ is a local minimizer of $E$, then $x^*$ is
also a local minimizer of $g-h$.
\end{prop}

The proof is given in Appendix \ref{append B}.

Based on the well-known results of d.c. optimization theory, if $x^*$ is a local
minimizer of $g-h$, then it must satisfy $\partial h(x^*)\subseteq\partial g(x^*)$(e.g.,
see\cite{HiriartUrruty1989,dualDCA1988,toland1979}), i.e.,
$E(x^*)\partial B(x^*)\subseteq\partial T(x^*)$. We emphasize that the function
$h$ depends on a given point $x^*$, which is not reflected in the function
symbols to keep the notation light.
Thus $x^*$ is a critical point of Problem $P_1 (P_2)$ according to the
following definition.
\begin{de}\label{critical point}
Let $x^*\in\mathcal{S}\cap\Omega$, $g$ and $h$ are defined as \eqref{define-gh}.
We say that $x^*$ is a critical point of Problem $P_1 (P_2)$ if
$\partial g(x^*)\cap\partial h(x^*)\neq\emptyset$, i.e.,
$E(x^*)\partial B(x^*)\cap\partial T(x^*)\neq\emptyset$.
If there exists $\epsilon>0$
such that $x^*$ is the only critical point in $\mathcal{B}_
{\scriptscriptstyle\mathcal{S}\cap\Omega}(x^*,\epsilon)$,
$x^*$ is said to be an isolated critical point of Problem $P_1(P_2)$.
\end{de}
We see that if a local minimizer of $E$ remains an isolated critical point
of Problem $P_1(P_2)$, it must be a strict local minimizer.
Upon analyzing the above, it is evident that any local minimizer of $E$ must be a
critical point of Problem $P_1(P_2)$.
We next present a sufficient condition for the local minimizer of $E$.
\begin{thm}
Let $x^*\in\mathcal{S}\cap\Omega$, $g$ and $h$ are defined as \eqref{define-gh}.
If there exists an $\epsilon>0$ such that $\partial h(x)\cap\partial g(x^*)\neq
\emptyset,\forall x\in\mathcal{B}_{\scriptscriptstyle\mathcal{S}\cap\Omega}(x^*,\epsilon)$,
then $x^*$ is a local minimizer of $E$.
\end{thm}
\prf The proof follows the arguments in \cite{Pham1998}.
For any $x\in\mathcal{B}_{\scriptscriptstyle\mathcal{S}\cap\Omega}(x^*,\epsilon)$,
according to the assumption,
there exists $v\in\partial h(x)\cap\partial g(x^*)$. $v\in\partial h(x)$ implies that
\begin{equation}\label{conj-h}
h(x)+h^*(v)=x^\top v\leq g(x)+g^*(v).
\end{equation}
Similarly, $v\in\partial g(x^*)$ implies that
\begin{equation}\label{conj-g}
g(x^*)+g^*(v)=(x^*)^\top v\leq h(x^*)+h^*(v).
\end{equation}
Combing \eqref{conj-h} and \eqref{conj-g}, we have
$$g(x)-h(x)\geq h^*(v)-g^*(v)\geq g(x^*)-h(x^*)=0,$$
namely, $T(x)-E(x^*)B(x)\geq0$. Hence, $E(x)\geq E(x^*), \forall x\in\mathcal{B}_
{\scriptscriptstyle\mathcal{S}\cap\Omega}(x^*,\epsilon)$.
\bbox

It's easy to see $x^*\in\Omega\cap\mathcal{S}$ is a solution to Problem $P_1(P_2)$
if and only if $g(x)-h(x)\geq0, \forall x\in\br^n$, which implies that $0$ is a
global minimizer of $g-h$, and therefore a local minimizer of $g-h$, owing
to the positive homogeneity of both
$g$ and $h$. Hence, it follows that $\partial h(0)\subseteq\partial g(0)$, i.e.,
\begin{equation}\label{suffcond}
E(x^*)\partial B(0)\subseteq\partial T(0).
\end{equation}
The following theorem indicates that \eqref{suffcond} is also
a sufficient condition for the solution of Problem $P_1(P_2)$.
\begin{thm}\label{equ cond P3}
Let $x^*\in\Omega\cap\mathcal{S}$, $g,h$ are defined as \eqref{define-gh}. The following statements are equivalent to each other:

(1) $E(x^*)\partial B(0)\subseteq\partial T(0)$;

(2) $x^*$ is a global minimizer of $E$, i.e., $E(y)\geq E(x^*),\forall y\in\Omega$;

(3) $0$ is a local minimizer of $g-h$.
\end{thm}
\prf
For convenience, we adopt the circular proof method.

$(1)\Rightarrow(2):$ The subdifferential inclusion implies that $E(x^*)v\in\partial T(0),~\forall v\in\partial B(0)$, by the definition of $\partial T$, for any $v\in\partial B(0)$, we have
$$T(y)\ge E(x^*)v^\top y,~\forall y\in\br^n.$$
Dividing both sides by $B(y)$ we obtain
$$E(y)\ge E(x^*)v^\top\frac y{B(y)},~\forall y\in\Omega.$$
Therefore, for any $y\in\Omega$, it also holds
\begin{equation}\label{subdifcontain}
E(y)\ge E(x^*)v^\top\frac y{B(y)},~\forall v\in\partial B(0).
\end{equation}
By the definition of $\partial B(0)$, for any $y\in\Omega$, we have
$$v^\top\frac y{B(y)}\leq1,~\forall v\in\partial B(0).$$
Moreover, the equality can be attained for all $y\in\Omega$. In fact,
according to the differential theory of convex function
(R.T.Rockafellar \cite[Th 23.1,Th 23.4]{1997Convex}), for any $y\in\br^n$, the one-sided directional derivative of $B$ at $0$ with respect to $y$ exists, denoted by $B'(0;y)$, and
$$\sup_{v\in\partial B(0)}v^\top y=B'(0;y)=\inf_{\lambda>0}\frac{B(\lambda y)-B(0)}{\lambda}=\inf_{\lambda>0}B (y)=B(y).$$
Therefore if $y\in\Omega$, it holds
$$\sup_{v\in\partial B(0)}\frac{v^\top y}{B(y)}=1.$$
As $\partial B(0)$ is a nonempty closed bounded set, there exists a corresponding $v_y\in\partial B(0)$
such that $v_{y}^\top\frac y{B(y)}=1$. Finally, by using inequality \eqref{subdifcontain}, we get the conclusion
$$E(y)\ge E(x^*)v_{y}^\top\frac y{B(y)}=E(x^*),~\forall y\in\Omega$$
as desired.

$(2)\Rightarrow(3):$ It is clear that $x^*$ minimizes $E$ if and only if
$0$ is a local minimizer of $g-h$.

$(3)\Rightarrow(1):$ According to the well known result of d.c. optimization,
if $0$ is a local minimizer of $g-h$, it holds
$\partial h(0)\subseteq\partial g(0)$\cite{HiriartUrruty1989,dualDCA1988,toland1979},
therefore the statement of $(1)$ is true.
\bbox

As a result, \eqref{suffcond} is an equivalent condition for the solution of
Problem $P_1 (P_2)$. Recall that the complete DCA introduced in Section
\ref{prelim} can be used to
reduce the value of the objective function $f_1-f_2$, and
if there is no strict decrease of the objective function we obtain $x^k$ that
satisfies $\partial f_2(x^k)\subseteq\partial f_1(x^k)$. Therefore,
given a candidate point $x^*\in\Omega$, if the complete DCA is conducted
one step starting from $x^0=0$ for $g-h$ defined by \eqref{define-gh}
and obtain $x^1$, then we must have $(g-h)(x^1)\leq(g-h)(x^0)=0$.
If $(g-h)(x^1)<0$, we obtain that the point $x^1$ satisfies $E(x^1)<E(x^*)$.
Otherwise, if $(g-h)(x^1)=0$, then $\partial h(0)\subseteq\partial g(0)$,
and consequently, according to Theorem \ref{equ cond P3}, we can conclude that
$x^*$ is a global minimizer of $E$.
Conversely, if $x^*$ is a global minimizer of $E$, i.e.,
$(g-h)(x)\geq0,\forall x\in\Omega$, it follows that $(g-h)(x^1)\geq0$.
On the other hand, owing to the
monotonicity of DCA, $(g-h)(x^1)\leq (g-h)(x^0)=0$. Therefore, $(g-h)(x^1)=0$.
To conclude, $(g-h)(x^1)=0$ if and only if $x^*$ is a global minimizer of $E$.
Unfortunately, the complete DCA needs to solve a convex maximization subproblem
in each iteration, which is difficult in general and time-consuming to perform.

In practice, one usually uses the simplified form of DCA, which is described in
\eqref{generDCA}. Although the simplified DCA still ensures that the objective
function monotonically decrease, if
$(g-h)(x^{1})=(g-h)(x^0)=0$, it does not imply that
$\partial h(0)\subseteq\partial g(0)$ but rather a weaker
one $\partial h(0)\cap\partial g(0)\neq\emptyset$. Therefore, if the
simplified DCA terminates, $x^*$ is not guaranteed to be a minimizer of $E$.

\section{The proximal-subgradient-d.c. algorithm (PS-DCA) for solving Problem
$P_1(P_2)$}
\label{alg-PS-DCA}
In this section, we devise an iterative algorithm
to obtain the critical point of Problem $P_1(P_2)$,
according to the global optimality condition.
For simplicity's sake, we introduce the notion of proximity operators. Assume $f\in\Gamma(\br^n)$,
$\prox_{f}:\br^n\rightarrow\br^n$ is the proximal operator of
$f$, which is defined as
\begin{equation}\label{eq:prox=argmin}
\prox_{f}(x):=\argmin_{y\in\br^n}
\Big(f(y)+\frac{1}{2}\|x-y\|^2\Big).
\end{equation}
Below, we present an equivalent condition concerning the critical points
of Problem $P_1(P_2)$.
\begin{thm}\label{th:important-p-2}
Let $x^*\in\mathcal{S}\cap\Omega$, $\lambda>0$. Then $x^*$ is a critical point
of Problem $P_1(P_2)$ if and only if there exists $v^*\in\partial B(x^*)$ such that $x^*=\prox_{\lambda T}(x^*+\lambda E(x^*)v^*)$.
\end{thm}
\prf
According to Definition \ref{critical point}, $x^*$ is a critical point of Problem
$P_1(P_2)$ if and only if there exists $v^*\in\partial B(x^*)$ such that $E(x^*)v^*\in\partial T(x^*)$, i.e.,
$\frac{1}{\lambda}\big[(x^*+\lambda E(x^*)v^*)-x^*\big]\in\partial T(x^*)$, which is equivalent to
$$
x^*=\prox_{\lambda T}(x^*+\lambda E(x^*)v^*).
$$
\bbox

Theorem \ref{th:important-p-2} suggests the iterative scheme
\begin{equation}\label{pgsa prox}
x^{k+1}=\prox_{\lambda_kT}(x^k+\lambda_kE(x^k)v^k),
\end{equation}
where $v^k\in\partial B(x^k)$. It is easy to see that \eqref{pgsa prox}
coincides with \eqref{pgsa}.
Given $x^k\in\Omega$, set the value of parameter $\lambda_k$. Let
\begin{equation}\label{def-til-gh}
{g}^k(x)=T(x)+\frac1{2\lambda_k}\|x\|^2,~~{h}^k(x)
=E(x^k)B(x)+\frac1{2\lambda_k}\|x\|^2.
\end{equation}
Then \eqref{pgsa prox} can be viewed as the simplified DCA applied to
the successive d.c. decomposition ${g}^k-{h}^k$. To see this, setting $u^0=x^k$
and running the simplified DCA one step, we obtain
$$ y^0\in\partial{h}^k(x^k), u^{1}\in\partial({g}^k)^*(y^0).$$
Since $\ri(\dom h)=\mathrm{\bf int}(\dom g)=\br^n$, it follows that $({g}^k)^*$
is differentiable on $\br^n$ and $\partial{h}^k(x)=E(x^k)\partial B(x)+\frac
x{\lambda_k}, (\forall x\in\br^n)$\cite[Th23.8]{1997Convex}. Therefore $y^0
\in\partial{h}^k(x^k)$ implies that there exists $v^k\in\partial B(x^k)$ such that $y^0=E(x^k)v^k+\frac{x^k}{\lambda_k}$. On the other hand, $u^1\in\partial({g}^k)^*(y^0)\Leftrightarrow y^0\in\partial{g}^k(u^{1})=\partial T(u^{1})+\frac{u^{1}}{\lambda_k} $. Therefore, $E(x^k)v^k+\frac{x^k}{\lambda_k}\in \partial T(u^{1})+\frac{u^{1}}{\lambda_k}$, i.e., $\lambda_k E(x^k)v^k+x^k\in(\lambda_k\partial T+I)(u^{1})$, this is equivalent to
$$
u^{1}=\prox_{\lambda_kT}(x^k+\lambda_kE(x^k)v^k).
$$
Hence, \eqref{pgsa prox} is actually a version of DCA.
Since the value of the regularization
parameter $\lambda_k$ in \eqref{pgsa prox} is related to the d.c. decomposition
of $T(x)-E(x^k)B(x)$, it may vary as $x^k$ changes,
which is different from PGSA\cite{PGSA2022}.

Our proposed method comprises two primary steps:
first, computing a candidate point, and refining it subsequently. More specifically,
in each iteration, \eqref{pgsa prox} is used to find a
candidate point, denoted by $l^k$, and let
\begin{equation}\label{param-rho}
\tilde{g}^k(x)=T(x)+\frac{\rho_k}{2}\|x\|^2,~~\tilde{h}^k(x)
=E(l^k)B(x)+\frac{\rho_k}{2}\|x\|^2,
\end{equation}
where $\rho_k>0$ is the regularization parameter.
Then a one-step simplified DCA is applied on $\tilde{g}^k-\tilde{h}^k$
starting from $x^0=0$ to possibly improve the candidate point $l^k$, i.e., we set
\begin{equation}\label{onestepdca}
s^k\in\partial\tilde{h}^k(0), t^{k}\in\partial(\tilde{g}^k)^*(s^k),
\end{equation}
then $t^k$ satisfies $(\tilde{g}^k-\tilde{h}^k)(t^k)\leq0$. If
$(\tilde{g}^k-\tilde{h}^k)(t^k)<0$, i.e., $E(t^k)<E(l^k)$, which means that
the DCA finds a better point, then the candidate point will be replaced,
otherwise it will be reserved. Due to the strong convexity of $\tilde{g}^k$
and $\tilde{h}^k$, $(\tilde{g}^k-\tilde{h}^k)(t^k)=0$
if and only if $t^k=0$. Since $s^k\in\partial\tilde{h}^k(0)$,
there exists $w^k\in\partial B(0)$ such that $s^k=E(l^k)w^k$.
Therefore $t^k\in\partial(\tilde{g}^k)^*(s^k)$ implies that
$E(l^k)w^k\in\partial T(t^k)+\rho_k t^k$, from which we deduce that
$t^k=\prox_{\frac{T}{\rho_k}}(\frac{E(l^k)}{\rho_k}w^k)$.
Since \eqref{pgsa prox} involves the proximal operator of $T$,
the subgradient of $B$, and the proposed algorithm incorporates DCA,
we refer to it as proximal-subgradient-d.c. algorithm (PS-DCA),
which is described in Algorithm \ref{alg:pg-dca}.

Now we make some remarks on Algorithm \ref{alg:pg-dca}.
It is known that evaluating the proximal operator of a function
$f\in\Gamma(\br^n)$ is a convex optimization problem,
and one can choose suitable algorithms and off-the-shelf solvers according to the
properties of $f$ (see \cite[Chapter 6]{Parikh2014Proximal} for more details of
these methods). The tolerance parameter $\epsilon'>0$ controls when to
accept $t^k$ found by the DCA.
It is worth to mention that as the d.c. components \eqref{param-rho}
are strongly convex, the only difference between the complete DCA and the simplified DCA lies on the choice of $s^k\in\partial\tilde{h}(0)$.
Therefore if $l^k$ is not a minimizer of $E$, whether it can be improved mainly
depends on the choice of $w^k\in\partial B(0)$. We will further explore how to select an appropriate $w^k$ to possibly improve
$l^k$ according to $T(x)$ and $B(x)$ in Section \ref{apply to GFMs prob}.
\begin{algorithm}[H]
\caption{PS-DCA:Proximal-subgradient-d.c. algorithm for solving $P_1(P_2)$}
\label{alg:pg-dca}
\begin{algorithmic}[1]
\State Input an initial point $x^0\in\mathcal{S}\cap\Omega$,
acceptance parameter $\epsilon'>0$.
\For{$k=0,1,\cdots$}
    \State Take $\lambda_k>0$;
    \State $v^k\in\partial B(x^k)$;
    \State $l^k=\prox_{\lambda_kT}(x^k+\lambda_k E(x^k)v^k)$;
    \State Take $\rho_k>0$;
    \State $w^k\in\partial B(0)$;
    \State $t^k=\prox_{\frac{T}{\rho_k}}(\frac{E(l^k)}{\rho_k}w^k)$;
    \State \textbf{if}~$(\tilde{g}^k-\tilde{h}^k)(t^k)<-\epsilon'$~
    \textbf{then} ~$y^k=t^k$;
    \State \textbf{else} ~ $y^k=l^k$;
    \State $x^{k+1}=\frac{y^k}{\|y^k\|}$;
    \State $k=k+1$;
    \EndFor
\end{algorithmic}
\end{algorithm}
\noindent

\section{Convergence analysis of PS-DCA}
\label{converg-of-PSDCA}
This section is dedicated to establishing convergence results of PS-DCA.
We begin by demonstrating the monotonicity and compactness for PS-DCA.
\begin{lem}[Monotonicity of PS-DCA]\label{th:decreasingofE(xk)}
Let $x^0\in\mathcal{S}\cap\Omega$ and $\epsilon'>0$. Suppose that the sequence
of iterates $(l^k,t^k,y^k,x^{k+1})$ are generated by Algorithm \ref{alg:pg-dca}.
Then for $k=0,1,2,\cdots$, $l^k\in\Omega,x^{k+1}\in\mathcal{S}\cap\Omega$ and if $y^k=t^k$, it holds
\begin{equation}\label{eq:th:decreasingofE(xk)}
E(x^k)\geq E(x^{k+1})+\frac{\rho_k\|t^k\|^2}{B(t^k)}  +\frac{1}{\lambda_kB(l^k)}\|x^k-l^k\|^2.
\end{equation}
Otherwise $y^k=l^k$, we have
\begin{equation}\label{decreaofsd}
E(x^k)\geq E(x^{k+1})+\frac{1}{\lambda_kB(l^k)}\|x^k-l^k\|^2.
\end{equation}
Thus, $E(x^k)=E(x^{k+1})$ if and only if $x^k=l^k=y^k=x^{k+1}$.
\end{lem}
\prf
Given that $x^k\in\mathcal{S}\cap\Omega$ and $\lambda_k>0$. Since $l^k$ can be
considered as the result of applying the simplified DCA to successive d.c.
decomposition ${g}^k-{h}^k$, we have, according to Theorem \ref{simp dca converg thm},
$$T(l^k)-E(x^k)B(l^k)=({g}^k-{h}^k)(l^k)\leq-\frac1{\lambda_k}
\|x^k-l^k\|^2.
$$
If $l^k\notin\Omega$, i.e., $B(l^k)=0$, we must have $T(l^k)=0$ and $x^k=l^k$, which is a contradiction since $x^k\in\Omega$, and thus $l^k\in\Omega$.
We divide $B(l^k)$ in the last inequality to obtain
\begin{equation}\label{E(xk),E(lk)}
E(x^k)\geq E(l^k)+\frac{1}{\lambda_kB(l^k)}\|x^k-l^k\|^2.
\end{equation}
If $y^k=l^k$, then
$x^{k+1}=\frac{l^k}{\|l^k\|}\in\mathcal{S}\cap\Omega$,
inequality \eqref{decreaofsd} follows from $E(l^k)=E(y^k)=E(x^{k+1})$. Otherwise
$y^k=t^k$, which implies that $(\tilde{g}^k-\tilde{h}^k)(t^k)<-\epsilon'$, consequently
$B(t^k)\neq0$ and $x^{k+1}=\frac{t^k}{\|t^k\|}\in\mathcal{S}\cap\Omega$.
Furthermore, according to Theorem \ref{simp dca converg thm}, it holds that
\begin{equation}\label{onedca}
T(t^k)-E(l^k)B(t^k)=(\tilde{g}^k-\tilde{h}^k)(t^k)\leq-\rho_k\|t^k\|^2.
\end{equation}
We divide $B(t^k)$ in \eqref{onedca} to obtain
$E(l^k)\geq E(t^k)+\frac{\rho_k\|t^k\|^2}{B(t^k)}$. Substituting this inequlity into \eqref{E(xk),E(lk)} and utilizing $E(t^k)=E(y^k)=E(x^{k+1})$, we obtain
\eqref{eq:th:decreasingofE(xk)}. The proof is complete.
\bbox
\begin{lem}[Compactness of PS-DCA]\label{th:compactness-new}
Let $x^0\in\mathcal{S}\cap\Omega$ and $\epsilon'>0$.
Suppose the sequence of
iterates $(l^k,t^k,y^k,x^{k+1})$ are generated by Algorithm \ref{alg:pg-dca}.
Assume the sequence $\{\lambda_k\}$ in Algorithm \ref{alg:pg-dca} has positive upper and lower bounds, $\{\rho_k\}$ has a positive lower bound. Then $\{l^k\}$ and $\{t^k\}$ are all bounded. Moreover, we have
$$
\|l^{k}-x^{k}\|\to0,~~~~\|x^{k}-x^{k+1}\|\to0.
$$
\end{lem}
\prf
Firstly, it is obvious that $\|x^k\|=1,~\forall k\in\bn$.
\par
As $v^k\in\partial B(x^k)$, we have
$$
B(z)\ge B(x^k)+(v^k)^\top(z-x^k),~~~~\forall z\in\br^n.
$$
Let $z=v^k+x^k$ to obtain
$B(v^k+x^k)\ge B(x^k)+(v^k)^\top v^k$, therefore,
$$\|v^k\|^2\le B(v^k+x^k)-B(x^k)\le B(v^k)$$
by using convexity and positive homogeneity of $B$.
If $\|v^k\|\neq0$, dividing by it in the last inequality yields
$\|v^k\|\le\beta$, where $\beta:=\max_{x\in\mathcal{S}}B(x)$.
\par
Since the proximal mapping $\prox_{f}$ is nonexpansive for $f\in\Gamma(\br^n)$ and $0=\prox_{\lambda_kT}(0)$, we obtain
$$\|l^k\|=\|\prox_{\lambda_kT}(x^k+\lambda_kE(x^k)v^k)-\prox_{\lambda_kT}(0)\|
\le\|x^k+\lambda_kE(x^k)v^k\|.$$
The boundedness of $\{l^k\}$ then follows from the facts that $\{v^k\}$, $\{\lambda_k\}$ are all bounded and the monotonicity of $\{E(x^k)\}$.
For the boundedness of $\{t^k\}$, note that the set $\partial B(0)$ is bounded,
then there exists a bounded set $S\subset\br^n$ so that
$\partial\tilde{h}(0)=E(l^k)\partial B(0)\subseteq S, (\forall k\in\mathbb{N})$.
Let $\partial\tilde{g}^*(S):=\bigcup\{
\partial\tilde{g}^*(s), s\in S\}$, then $\partial\tilde{g}^*(S)$ is a nonempty bounded set\cite[Th24.7]{1997Convex}, thus $\{t^k\}\subset\partial\tilde{g}^*(S)$ is bounded.
\par
For the second statement, as $\{E(x^k)\}$ is decreasing and nonnegative, it must
converge. According to Lemma \ref{th:decreasingofE(xk)},
$$\|x^{k}-l^{k}\|^2\le\lambda_{k}B(l^{k})[E(x^{k})-E(x^{k+1})]\to 0,~~~~k\to\infty,$$
which implies that $\|x^{k}-l^{k}\|\to 0, k\to\infty$.
The proof of the final conclusion is grounded on the fact pertaining to PS-DCA:
there exists a $K\in\mathbb{N}$ sufficiently large so that $y^k=l^k, \forall k\geq K$. Denote $\mathcal{K}:=\{k\in\mathbb{N}\mid
y^k=t^k\}$. It suffices to show that $\mathcal{K}$ is a finite set. If not,
letting $k\in\mathcal{K}$ tend to infinity, we have
\begin{equation}\label{tk-mono}
\frac{\lambda_{k}\rho_{k}B(l^k)}{B(\frac{t^{k}}{\|t^{k}\|})}\|t^{k}\|+\|x^{k}-l^{k}\|^2\leq
\lambda_{k}B(l^{k})[E(x^{k})-E(x^{k+1})]\to 0
\end{equation}
by using \eqref{eq:th:decreasingofE(xk)}. Note that $\{B(l^k)\}$ has a positive lower bound, i.e.,
there exists a constant $c>0$ such that
$B(l^k)>c,(\forall k\in\mathbb{N})$.
Otherwise, there must exist a subsequence of $\{B(l^k)\}$ that tends towards $0$.
For simplicity, we assume $B(l^{k})\rightarrow0,(k\rightarrow\infty)$.
By using the monotonicity of PS-DCA, we have
$E(l^k)=\frac{T(l^k)}{B(l^k)}\leq E(x^0),(\forall k\in\mathbb{N})$,
thus $T(l^{k})\rightarrow0,(k\rightarrow\infty)$.
Given that $\{l^k\}$ is bounded, it must have a convergent subsequence, denoted
by ${l^{k_j}}$, and assume
that $l^{k_j}\rightarrow a,(j\rightarrow\infty)$.
As a result, $B(l^{k_j})\rightarrow B(a)$
and $T(l^{k_j})\rightarrow T(a)$, from which we deduce that $T(a)=0$ and $B(a)=0$.
This implies that $a=0$ according to our basic assumptions on $T$ and $B$.
Since any convergent subsequence of $\{l^k\}$
converges to $0$, it must have $l^k\rightarrow0$, which contradicts the facts
that $\|l^k-x^k\|\rightarrow0$ and $\|x^k\|=1$. Therefore, we derive
$t^{k}\to 0,(k\in\mathcal{K},k\rightarrow\infty)$ from \eqref{tk-mono},
providing $\frac1{B(\frac{t^k}{\|t^k\|})}\geq\frac1\beta$ with $\beta>0$
and the assumptions for sequences $\{\lambda_k\}$ and
$\{\rho_k\}$. This fact contradicts
$(\tilde{g}^k-\tilde{h}^k)(t^k)<-\epsilon',(\forall k\in\mathcal{K})$
for a given $\epsilon'>0$. Therefore the set $\mathcal{K}$ contains at most
a finite number of elements and the inequality
$$\|x^{k}-x^{{k}+1}\|\leq
\|x^{k}-l^{k}\|+\big\|l^{k}-\frac{l^{k}}{\|l^{k}\|}\big\|\to 0,~~~~k\to\infty$$
then follows from the fact that $\|l^{k}\|\to \|x^{k}\|=1,(k\to\infty)$.
\bbox

Based on the monotonicity and compactness of Algorithm \ref{alg:pg-dca},
we are ready to establish its subsequential convergence property.
\begin{thm}\label{th:global-convergence}
Suppose that $\{\lambda_k\}$ in Algorithm \ref{alg:pg-dca} has positive upper and lower bounds, $\{\rho_k\}$ has a positive lower bound. Then

(1). Any accumulation point of the sequence $\{x^k\}$ generated by Algorithm
\ref{alg:pg-dca} is a critical point of Problem $P_1(P_2)$.

(2). Either the sequence $\{x^k\}$ converges or the set of accumulation points
of $\{x^k\}$ form a continuum on $\mathcal{S}\cap\Omega$,
i.e., the set of accumulation points possesses the cardinality of a continuum,
and it does not contain any isolated accumulation points.
\end{thm}
\prf
(1). As $\{x^k\}$, $\{v^k\}$ and $\{\lambda_k\}$ are all bounded, there must exist a convergent subsequence for all of the sequences, we may assume that
$$x^k\to x^*,~~~v^k\to v^*,~~~\lambda_k\to\lambda^*,~~~~as~k\to\infty.$$
Lemma \ref{th:compactness-new} suggests that $l^k\to x^*(k\to\infty)$.
It is easy to see that $x^*\in\mathcal{S}$, we next show $x^*\in\Omega$.
According to Lemma \ref{th:decreasingofE(xk)}, we have
\begin{equation}\label{Mono E}
E(x^k)\geq E(x^{k+1})+\frac{1}{\lambda_kB(l^k)}\|x^k-l^k\|^2,~\forall k,
\end{equation}
from which we deduce that
\begin{equation}\label{xk-lk}
\frac{1}{\lambda_kB(l^k)}\|x^k-l^k\|^2\to0,(k\to\infty).
\end{equation}
On the other hand, from \eqref{Mono E}, we also have
$$T(x^{k+1})\leq E(x^k)B(x^{k+1})
-\frac{B(x^{k+1})}{\lambda_kB(l^k)}\|x^k-l^k\|^2,~\forall k.$$
By letting $k\to\infty$ and combining with \eqref{xk-lk}, we obtain
$T(x^*)\leq c^*B(x^*)$, where $c^*=\lim_{k\to\infty}E(x^k)$. Therefore,
if $B(x^*)=0$, then $T(x^*)=0$ and thus by assumption we deduce $x^*=0$,
which is a contradiction. Consequently, $x^*\in\Omega$ and $c^*=E(x^*)$ by
the continuity of $E$.

Denote
\begin{equation}\label{f def}
f(l,\lambda, x, v):=T(l)+\frac{1}{2\lambda}\|l-(x+\lambda E(x)v)\|^2,
\end{equation}
$l^k=\prox_{\lambda_kT}(x^k+\lambda_k E(x^k)v^k)$ implies
$f(l^k,\lambda_k,x^k,v^k)\le f(l, \lambda_k,x^k,v^k)$ for any $l\in\br^n$. Let $k\to\infty$ to obtain
$f(x^*,\lambda^*,x^*,v^*)\le f(l, \lambda^*,x^*,v^*)$, i.e.,
$x^*=\prox_{\lambda^*T}(x^*+\lambda^* E(x^*)v^*)$. Note that $v^k\in\partial B(x^k)$ and the graph of $\partial B$ is closed on $\br^n\times\br^n$\cite[Th24.4]{1997Convex}, therefore
$v^*\in\partial B(x^*)$ and thus $x^*$ is a critical point of Problem $P_1(P_2)$
according to Theorem \ref{th:important-p-2}.

(2). Since the sequence $\{x^{k}\}$ generated by Algorithm \ref{alg:pg-dca}
is bounded and $\|x^{k}-x^{k+1}\|\to0,(k\to\infty)$ according to Lemma
\ref{th:compactness-new}, then $\{x^k\}$ converges, or its accumulation points
form a continuum (see \cite[Theorem 26.1]{globalconverg}).
\bbox

If the one-step DCA starting from $0$ is removed, PS-DCA reduces
to the proximal-subgradient-algorithm (PSA), which is stated in Algorithm
\ref{alg:psa}. The iteration scheme of PSA coincides with PGSA \eqref{pgsa}
except a projection onto the sphere is added in PSA to guarantee the boundedness
of the iterates. Clearly, PSA possesses monotonicity and compactness as well,
thus, the subsequential convergence property, as stated in Theorem
\ref{th:global-convergence}, also applies to PSA.
Apart from these properties, we also have the following results of PSA regarding
the local minimizers of $E$.
\begin{lem}\label{CLM property}
Suppose $x^*$ is a local minimizer of $E$,
$\{x^k\}\subset\mathcal{S}\cap\Omega$ is any sequence satisfying
$x^k\rightarrow x^*$. For each $x^k$, define $v^k,l^k$ according to PSA,
denote $z^k:=l^k/\|l^k\|$. If $\{\lambda_k\}$ has positive upper and lower bounds, then $z^k\rightarrow x^*$.
\end{lem}
\prf
The proof is inspired by \cite[Lemma 7]{NIPS2012}.
The boundedness of $\{x^k\},\{\lambda_k\}$ and $\{v^k\}$ yield $\{l^k\}$ is
bounded, therefore, any subsequence of $\{l^k\}$ has a further convergent
subsequence $l^{k_i}\rightarrow l^*$. Moreover, we may assume that
$\lambda_{k_i}\rightarrow\lambda^*$, $v^{k_i}\rightarrow v^*$. Then $\lambda^*>0$ and $v^*\in\partial B(x^*)$.
By the definition of $l^k$, we have $f(l^{k_i},\lambda_{k_i},x^{k_i},v^{k_i})\leq f(l,\lambda_{k_i},x^{k_i},v^{k_i}),\forall l\in\br^n$, where $f$ is defined as \eqref{f def}. Let $i\rightarrow\infty$ to obtain
$l^*=\prox_{\lambda^*T}(x^*+\lambda^* E(x^*)v^*)$. Note that $x^*$ is a local minimum of $E$, it must satisfy $E(x^*)\partial B(x^*)\subseteq\partial T(x^*)$, i.e., for any $v^*\in\partial B(x^*)$, it holds
$0\in\partial T(x^*)-E(x^*)v^*=\partial T(x^*)+\frac1{\lambda^*}
(x^*-x^*-\lambda^* E(x^*)v^*)$, which is equivalent to $x^*=
\prox_{\lambda^*T}(x^*+\lambda^* E(x^*)v^*)$. Therefore $l^*=x^*$.
Since any subsequence of $\{l^k\}$ has a further subsequence that converges to
$x^*$, the entire sequence converges to $x^*$. Consequently, $z^k\rightarrow x^*$.
\bbox

Based on Lemma \ref{CLM property}, along with the monotonicity and subsequential
convergence of PSA, we can establish the closedness
property of PSA concerning local minimizers of $E$, which is presented as
follows:
\begin{thm}[Closedness of PSA]\label{local-converge}
Let $x^*$ denote a local minimizer of $E$ that is further an isolated critical
point of Problem $P_1(P_2)$. Suppose $\{x^k\}$ is generated by PSA. If
$\{\lambda_k\}$ has positive upper and lower bounds, then
$x^0\in\mathcal{B}_\delta(x^*)$ with $\delta>0$ sufficiently small implies
$x^k\rightarrow x^*$.
\end{thm}

Since Theorem \ref{local-converge} can be demonstrated using techniques
similar to those employed in \cite[Lemma 10]{NIPS2012} and
\cite[Theorem 1]{NIPS2012}, we omit the detailed proof here.

Theorem \ref{local-converge} indicates that as long as
an iterate is sufficiently near a local minimizer of $E$, the entire sequence
generated by PSA converges to this local minimizer. However, PS-DCA generally lacks closure at local minimizers of
$E$, which can be seen as an advantage to a certain degree. In fact,
if an iterate lies near a local minimizer of $E$, and the objective value at
this local minimizer is relatively high, the sequence generated by PSA is
likely to become trapped in the neighborhood of this local minimizer due to the
algorithm's enclosure of local minimizers. In contrast,
a one-step DCA initiated from zero allows subsequent iterates of PS-DCA to escape
from this neighborhood and move towards a higher-quality local minimizer.
As observed in the experiments of Section \ref{exp}, the sequence generated by
PS-DCA is not drawn to the local minimizers with high objective values.

\begin{algorithm}[H]
\caption{PSA:Proximal-subgradient algorithm for solving Problem $P_1(P_2)$}
\label{alg:psa}
\begin{algorithmic}[1]
\State Input an initial point $x^0\in\mathcal{S}\cap\Omega$.
\For{$k=0,1,\cdots$}
    \State Take $\lambda_k>0$;
    \State $v^k\in\partial B(x^k)$;
    \State $l^k=\prox_{\lambda_kT}(x^k+\lambda_k E(x^k)v^k)$;
    \State $x^{k+1}=\frac{l^k}{\|l^k\|}$;
    \State $k=k+1$;
    \EndFor
\end{algorithmic}
\end{algorithm}
\noindent

Finally, we demonstrate that in practice, PS-DCA or PSA is able to yield
an approximate critical point for Problem $P_1(P_2)$.
When the algorithm stops, that is
$E(x^k)-E(x^{k+1})\leq\epsilon^2$ with $\epsilon>0$ sufficiently small,
according to Lemma \ref{th:decreasingofE(xk)},
we have $\|l^k-x^k\|^2\leq\lambda_kB(l^k)[E(x^k)-E(x^{k+1})]\leq\epsilon^2
\lambda_kB(l^k)$. Under the assumptions of Lemma \ref{th:compactness-new},
$\{\lambda_k\}$ and $\{l^k\}$ are all bounded, thus there exists a constant
$c>0$ such that $\|l^k-x^k\|\leq c\epsilon$.
Since $l^k$ is actually obtained by performing the one-step simplified DCA on
$g^k-h^k$ starting from $x^k$,
where $g^k$ and $h^k$ are defined as \eqref{def-til-gh}, therefore
$l^k\in\partial({g}^k)^*(y^0)$ with $y^0\in\partial{h}^k(x^k)$, this implies that
\begin{equation}\label{lk-inte-xk}
\partial{h}^k(x^k)\cap\partial{g}^k(l^k)\neq\emptyset.
\end{equation}
On the other hand, by using the continuity of the subdifferential mapping
$\partial g^k$, for any $\alpha>0$, there exists a corresponding $\delta>0$ such that
$\partial{{g}^k}(z)\subseteq\partial{{g}^k}(x^k)+\alpha\mathcal{B}(0,1),~\forall
\|z-x^k\|\leq\delta$\cite[Cor24.5.1]{1997Convex}. Therefore,
for any given $\alpha>0$, there exists a corresponding $\delta>0$ such that when
$\epsilon$ is sufficiently small, $c\epsilon\leq\delta$ holds true, then we have
\begin{equation}\label{continu-subgra}
\partial{{g}^k}(l^k)\subseteq\partial{{g}^k}(x^k)+\alpha\mathcal{B}(0,1).
\end{equation}
Combining \eqref{lk-inte-xk} and \eqref{continu-subgra}, we obtain
$\partial{{h}^k}(x^k)\cap(\partial{{g}^k}(x^k)+a\mathcal{B}(0,1))\neq\emptyset$.
Note that $x^k\in\Omega\cap\mathcal{S}$ is a critical point of Problem $P_1(P_2)$
if and only if $\partial{{h}^k}(x^k)\cap\partial{{g}^k}(x^k)\neq\emptyset$.
Therefore, $x^k$ can be regarded as an approximate critical point of Problem
$P_1(P_2)$.

\section{Application to the generalized graph Fourier mode problem}
\label{apply to GFMs prob}
In this section, we investigate global sequential convergence of the sequence
generated by PS-DCA or PSA for the generalized graph Fourier mode problem
aimed at minimizing the graph directed variation. Then we deduce the method for
selecting $w^k$ within PS-DCA for this type of problem.

\subsection{Global sequential convergence of PS-DCA or PSA}
Firstly, we review the definition of generalized graph Fourier modes
based on graph directed variation minimization.
For the sake of simplicity, we introduce the notation
$\bfone=(1,\cdots,1)^\top\in\br^n$.

\begin{de}\cite{2018irregular}
\label{de-GFM-gdv}
Assume $W=[w_{ij}]\in\br^{n\times n}$ is a weight matrix of a directed graph
that has no isolated point. Let $T$ be the graph directed variation operator,
i.e., $T(x)=\sum_{i,j=1}^nw_{ij}[x_i-x_j]_+$, where $[x]_+=\max(0,x)$.
$Q$ is a diagonal matrix of order $n$, with all its diagonal elements being
positive.
Then the set of $(T,Q)$-graph Fourier modes ($(T,Q)$-GFMs)
$\{u_k\}^{n}_{k=2}$ is defined as
\begin{equation}\label{de-T-Q-GFM}
\begin{split}
u_k:=\argmin_{x\in X}&~T(x)\\
        \text{s.t.}~~ &\|Q^{\frac12}x\|=1,
\end{split}
\end{equation}
where $X=\{x\in\br^n\mid U^\top_{k-1}Qx=0\}$,
$U_{k-1}=(u_1,\cdots,u_{k-1})$ with $u_1=\frac1{\|Q^{\frac12}\bfone\|}
\bfone$.
\end{de}

It is well known that a set of graph Fourier modes, derived by minimizing the
graph directed variation, is sparser than the Laplacian basis and is highly
effective in identifying clusters within the graph\cite{dirGFT2017,l1GFT2021}.
When $Q$ is taken to be the identity matrix $I$, the set of $(T,I)$-GFMs
minimizes the graph directed total variation as well,
thereby forming the graph Fourier basis for directed graph. Furthermore,
when the graph is undirected, $T(x)$ equals to $l_1$ norm variation and the set
of $(T,I)$-GFMs is precisely the $l_1$ Fourier basis.
It can be verified that for any $k$, \eqref{de-T-Q-GFM} is a
special case of problem $P_1$, and thus equivalent to the following fractional
minimization problem
\begin{equation}\label{GFM-fracprogram}
\min_{x\in\br^m\setminus\{0\}}\frac{T(Vx)}{\|Q^{\frac12}Vx\|},
\end{equation}
where $m=n-k$, $V=(v_1,\cdots,v_m)$ with $\{v_i\}^m_{i=1}$ is an
orthonormal basis for $X$. For the proof of
equivalence between problem \eqref{de-T-Q-GFM} and problem \eqref{GFM-fracprogram},
please refer to Appendix \ref{appendix}. Next, we present the main result of
this section, from which the global sequential convergence of the sequence
generated by PS-DCA or PSA for problem \eqref{GFM-fracprogram} can be immediately
derived according to Theorem \ref{th:global-convergence}.
Let
$\mathcal{N}:=\{x\in\mathcal{S}\mid T(Vx)=0\}$.
\begin{thm}\label{L1norm}
For the $(T,Q)$-GFMs defined as above. Given $1\leq k\leq n-1$, if
$\mathcal{N}$ is a finite set, then problem \eqref{GFM-fracprogram}
has only a finite number of critical points.
\end{thm}

\prf
Since $\|\cdot\|$ is differentiable on $\br^m\setminus\{0\}$,
we have $\partial\|\cdot\|(x)=\{x/\|x\|\},\forall x\neq0$.
Hence for $\forall x^*\in\mathcal{S}$, $v^*\in\partial\|Q^{\frac12}
V\cdot\|(x^*)=\{\frac{V^\top QVx^*}{\|Q^{\frac12}Vx^*\|}\}$\cite[Th 23.9]{1997Convex}.
Note that $x^*\in\mathcal{S}$ is a critical point of problem \eqref{GFM-fracprogram}
if and only if $\frac{T(Vx^*)}{\|Q^{\frac12}Vx^*\|}
v^*\in \partial T(V\cdot)(x^*)$, i.e.,
\begin{equation}\label{gdv-Q-subgrad}
T(Vx)\ge T(Vx^*)+\frac{T(Vx^*)}{\|Q^{\frac12}Vx^*\|^2}
(Vx^*)^\top QV(x-x^*)=\frac{T(Vx^*)}{\|Q^{\frac12}Vx^*\|^2}
(Vx^*)^\top QVx,~~~~\forall x\in\br^m.
\end{equation}
It is clear that if $x^*\in\mathcal{N}$, \eqref{gdv-Q-subgrad} always holds true.
Consequently, all elements in set $\mathcal{N}$ constitute critical points of
problem \eqref{GFM-fracprogram}.
Based on the assumption of the theorem, it is sufficient to demonstrate that
there exist only a finite number of critical points of problem \eqref{GFM-fracprogram}
within $\mathcal{S}\setminus\mathcal{N}$. In the following proof, we assume that
$T(Vx^*)>0$.

Let $Vx=y$, $Vx^*=y^*$,
then \eqref{gdv-Q-subgrad} becomes
\begin{equation}\label{gdv-Q-crtical}
T(y)\geq\frac{T(y^*)}{\|Q^{\frac12}y^*\|^2}
(y^*)^\top Qy,~~~~\forall y\in X,
\end{equation}
where $y^*\in\mathcal{S}\cap X$ and satisfies $T(y^*)>0$. Let $x_1=Q^{\frac12}y$,
$x^*_1=Q^{\frac12}y^*$. It can be seen that \eqref{gdv-Q-crtical} is equivalent to
\begin{equation}\label{gdv-Q-x1}
T(Q^{-\frac12}x_1)\geq T(Q^{-\frac12}x^*_1)\frac{(x^*_1)^\top x_1}{\|x^*_1\|^2},
~~~~\forall x_1\in X_1,
\end{equation}
where $X_1=\{x\in\br^n\mid U^\top_{k-1}Q^{\frac12}x=0\}$, $x_1^*\in X_1$
satisfying $(x_1^*)^\top Q^{-1}x_1^*=1$ and $T(Q^{-\frac12}x^*_1)>0$.
Let $V_1=(l_1,\cdots, l_m)$ be an orthonormal basis for $X_1$, assume $x_1=V_1z$
and $x_1^*=V_1z^*$. Therefore, \eqref{gdv-Q-x1} can be expressed as
$$T(Q^{-\frac12}V_1z)\geq T(Q^{-\frac12}V_1z^*)\frac{(z^*)^\top z}{\|z^*\|^2},
~~~~\forall z\in\br^m.$$
Let $A=Q^{-\frac12}V_1$, the above inequality is
\begin{equation}\label{gdv-Q-A}
T(Az)\geq T(Az^*)\frac{(z^*)^\top z}{\|z^*\|^2},
~~~~\forall z\in\br^m,
\end{equation}
where $z^*\in\br^m$ satisfies $(z^*)^\top V_1^\top Q^{-1}V_1z^*=1$ and
$T(Az^*)>0$.

Assume $A=\begin{bmatrix}
a_1^\top\\
\vdots\\
a_n^\top
\end{bmatrix}
$, where $a_i\in\br^m$. Then $Az=\begin{bmatrix}
a_1^\top z\\
\vdots\\
a_n^\top z
\end{bmatrix}$,
and
$$T(Az)=\sum_{i,j=1}^nw_{ij}[(a_i-a_j)^\top z]_+,~~~~\forall z\in\br^m.$$

When restricted in a subspace and a sufficiently small neighborhood of $z^*$,
$T(Az)$ can be expressed as a linear function of $z$.
To see this, based on the sign of $(a_i-a_j)^\top z^*$,
we introduce the following partition with respect to the set of indices set
$\Lambda=\{(i,j)\mid 1\leq i\leq n, 1\leq j\leq n\}$:
\begin{gather*}
 \Lambda_0:=\{(i,j)\in\Lambda\mid(a_i-a_j)^\top z^*=0\},\\
 \Lambda_+:=\{(i,j)\in\Lambda\mid(a_i-a_j)^\top z^*>0\},\\
 \Lambda_-:=\{(i,j)\in\Lambda\mid(a_i-a_j)^\top z^*<0\}.
\end{gather*}
It is easy to see that $\Lambda_+$ is nonempty as $T(Az^*)>0$. Denote
$\epsilon:=\min_{(i,j)\notin\Lambda_0}|(a_i-a_j)^\top z^*|$,
$\alpha:=\max_{(i,j)\notin\Lambda_0}\|a_i-a_j\|$. Note that
$a_i-a_j\neq\bf{0}$, $\forall(i,j)\notin\Lambda_0$, hence $\alpha>0$.
Let $H:=\{z\in\br^m\mid (a_i-a_j)^\top z=0,\forall (i,j)\in\Lambda_0\}$.
Clearly, $H$ is a linear subspace of $\br^m$. If $\Lambda_0$ is the empty set,
then $H=\br^m$. Furthermore, $H\neq\{\bf0\}$ since $z^*\in H$.

For $\forall\Delta z\in H$ and $\|\Delta z\|<\epsilon/\alpha$,
we have
$$|(a_i-a_j)^\top\Delta z|\leq\|a_i-a_j\|\|\Delta z\|<\alpha\cdot\frac\epsilon
\alpha=\epsilon\leq|(a_i-a_j)^\top z^*|, ~~\forall(i,j)\notin\Lambda_0.$$
Let
$z=z^*+\Delta z$, then $z\in H$. From the above inequality,
we can deduce that
\begin{itemize}
  \item for $\forall(i,j)\in\Lambda_0$, we have $(a_i-a_j)^\top z=0$, thus
  $[(a_i-a_j)^\top z]_+=0$.
  \item for $\forall(i,j)\in\Lambda_+$, we have $(a_i-a_j)^\top z>0$, thus
  $[(a_i-a_j)^\top z]_+>0$.
  \item for $\forall(i,j)\in\Lambda_-$, we have $(a_i-a_j)^\top z<0$,
  thus $[(a_i-a_j)^\top z]_-=0$.
\end{itemize}
Consequently, for $\forall\Delta z\in H$ and $\|\Delta z\|<\epsilon/\alpha$,
we have
\begin{equation}\label{linear-equa-T-Az}
\begin{split}
T(Az)&=\left(\sum_{(i,j)\in\Lambda_0}+\sum_{(i,j)\in\Lambda_+}+
\sum_{(i,j)\in\Lambda_-}\right)w_{ij}[(a_i-a_j)^\top z]_+\\
&=\sum_{(i,j)\in\Lambda_+}
w_{ij}[(a_i-a_j)^\top z]_+=f^\top z,
\end{split}
\end{equation}
where $f:=\sum_{(i,j)\in\Lambda_+}w_{ij}(a_i-a_j)$. In particular,
setting $\Delta z=0$, we obtain $T(Az^*)=f^\top z^*$.
Note that since $W$ and $A$ are fixed in the problem and there are only a
limited number of ways to partition $\Lambda$, thus $f$ can only take on a
finite number of possibilities.

Substitute \eqref{linear-equa-T-Az} into
\eqref{gdv-Q-A} and note that $z=z^*+\Delta z$, we obtain
$$
f^\top(z^*+\Delta z)\geq\frac{f^\top z^*}{\|z^*\|^2}(z^*)^\top(z^*+\Delta z)
,~~~~\forall \Delta z\in H,\|\Delta z\|<\epsilon/\alpha.
$$
Assume that $\dim H=h$. Let $M\in\br^{m\times h}$ be the matrix formed by a set of
orthonormal bases of $H$. Then $z^*=Mc^*$, $\Delta z=M\Delta c$, where
$\|c^*\|=\|z^*\|$, $\|\Delta c\|=\|\Delta z\|$. Therefore,
$$f^\top M(c^*+\Delta c)\geq\frac{f^\top Mc^*}{\|c^*\|^2}(c^*)^\top(c^*+\Delta c),
~~~~\forall \Delta c\in\br^h, \|\Delta c\|<\epsilon/\alpha.$$
Let $g=M^\top f$, we obtain
$$g^\top(c^*+\Delta c)\geq\frac{g^\top c^*}{\|c^*\|^2}(c^*)^\top(c^*+\Delta c).$$
Hence,
$$
(g^\top-\frac{g^\top c^*(c^*)^\top}{\|c^*\|^2})\Delta c\geq0,~~~~\forall\Delta c
\in\br^h,\|\Delta c\|<\epsilon/\alpha.
$$
Given the arbitrariness of $\Delta c$, it can be inferred that
$g^\top=\frac{g^\top c^*(c^*)^\top}{\|c^*\|^2}$. Therefore,
$$\|g\|^2=g^\top g=\frac{g^\top c^*(c^*)^\top g}{\|c^*\|^2}=
\frac{(g^\top c^*)^2}{\|c^*\|^2},$$
which implies
\begin{equation}\label{linear-g-c}
g^\top c^*=\|g\|\|c^*\|.
\end{equation}
Note that $g^\top c^*=f^\top M
c^*=f^\top z^*=T(Az^*)>0$, thus we can only adopt the positive sign.
\eqref{linear-g-c} indicates that $c^*=ag$, where $a$ is a positive number.
Remember that $z^*$ satisfies $(z^*)^\top V_1^\top Q^{-1}V_1z^*=1$ and $z^*=Mc^*$,
then we have $(c^*)^\top M^\top V_1^\top Q^{-1}V_1Mc^*=1$. This implies that
$\|Q^{-\frac12}V_1Mc^*\|=1$ and therefore $a=\frac1{\|Q^{-\frac12}V_1Mg\|}$.
Hence,
$$z^*=Mc^*=\frac{Mg}{\|Q^{-\frac12}V_1Mg\|}=\frac{MM^\top f}{\|Q^{-\frac12}V_1M
M^\top f\|}.$$
The conclusion then follows from the fact that $f$ has a
limited number of possibilities only, so there are
only a finite number of $z^*$ with the above representation.
\bbox

\textbf{Remark 1.} For an undirected and connected graph, $T(x)=0$
if and only if $x$ is a constant signal. In this case, the condition of
Theorem \ref{L1norm} is always satisfied. In the context of disconnected
undirected and directed graphs, when $k$ is small, $\mathcal{N}$ may
have an infinite number of elements, forming a continuum.
For instance, in the case of an undirected graph with
more than three connected components,
$T(x)=0$ if and only if $x$ is a constant signal across each connected component.
When $X=\{x\in\br^n\mid{\bf1}^\top Qx=0\}$, there exist
infinitely many such piecewise constant signals within $X\cap\mathcal{S}$,
collectively forming a continuum.
For the case of directed graphs, let's consider a simple example.
Let $T(x)=w_{12}(x_1-x_2)_+ +w_{23}(x_2-x_3)_+ +w_{34}(x_3-x_4)_+$,
which represents the graph directed variation of a path.
Given $X=\{x\in\br^4\mid{\bf1}^\top x=0\}$, one can verify that set
$V^\top\{x\in X\cap\mathcal{S}\mid x_1=x_2<x_3<x_4\}$ is a subset of $\mathcal{N}$ and
obviously, the elements within this set form a continuum.

\textbf{Remark 2.}
If the denominator in problem \eqref{GFM-fracprogram} is replaced with
$\ell_1$ norm, then the number of critical points in set
$\mathcal{S}\setminus\mathcal{N}$
is no longer finite. For example,
let $T(x)=|x_1-x_2|+|x_3-x_4|$, $B(x)=|x_1|+|x_2|+|x_3|+|x_4|$.
Then on set $\{x\in\br^4\mid x_1>0,x_2<0,x_3>0,x_4<0\}$, it holds that
$T(x)=B(x)>0$, thus $E(x)\equiv1$ and $\partial T(x)=\partial B(x)$. Therefore,
for any subspace $X$ in $\br^4$, all elements in set
$\{x\in X\cap\mathcal{S}\mid x_1>0,x_2<0,x_3>0,x_4<0\}$ satisfy
$\partial T(x)\cap E(x)\partial B(x)\neq\emptyset$,
making them critical points
of problem $\min_{x\neq0}\frac{T(x)}{B(x)}$
within $X$. However, provided that this set is nonempty, its elements can
always form a continuum.

Theorem \ref{L1norm} indicates that if the number of zeros of $T(x)$ within the
feasible region is finite, then the entire sequence generated by either PSA or
PS-DCA will converge to a critical point of the
generalized graph Fourier mode problem \eqref{de-T-Q-GFM}.
Since the quality of the solution returned by PS-DCA is affected by DCA's
capability to enhance the candidate point,
in the subsequent subsection, we develop a method for
selecting a suitable subgradient $w^k$ within PS-DCA to possibly improve
$l^k$ for this class of generalized graph Fourier mode problems.

\subsection{Choosing the subgradient $w^k$ in PS-DCA for problem \eqref{GFM-fracprogram}}
\label{choos-wk}
Denote by $d_{min}$ the smallest diagonal element of $Q$.
Let $E(x):=\frac{T(Vx)}{\|Q^{\frac12}Vx\|}$.
According to Algorithm \ref{alg:pg-dca}, $w^k\in\partial\|Q^{\frac12}V\cdot\|(0)=
V^\top Q^{\frac12}\partial\|\cdot\|(0)$,
thus choosing $w^k$ is equivalent to selecting $\tilde{v}_0\in\partial\|\cdot\|(0)$.
Let $\spa V$ denote the column space of $V$, i.e., $\{Vx\mid x\in\br^m\}$.
For convenience,
we choose $\tilde{v}_0$ such that $Q^{\frac12}\tilde{v}_0\in\spa V$. Note that
such $\tilde{v}_0$ does exists. In fact, let $\tilde{v}_0=Q^{-\frac12}Vc$,
it can be verified that whenever $\|c\|\leq\sqrt{d_{min}}$,
$\|\tilde{v}_0\|=\|Q^{-\frac12}Vc\|\leq
\|Q^{-\frac12}\|\|c\|\leq\frac1{\sqrt{d_{min}}}\sqrt{d_{min}}=1$,
which implies that $\tilde{v}_0\in\partial\|\cdot\|(0)$. With this choice of
$\tilde{v}_0$, $w^k=V^\top Q^{\frac12}\tilde{v}_0=V^\top Q^{\frac12}
Q^{-\frac12}Vc=c$. Choosing $w^k$ is thus equivalent to choosing $c$.

On the other hand, provided that $l^k$ obtained from the
proximal subgradient step in PS-DCA is not a minimizer of problem
\eqref{GFM-fracprogram}, according to Theorem \ref{equ cond P3}, there must exist
some $\tilde{w}^k\in E(l^k)\partial\|Q^{\frac12}V\cdot\|(0)$ such that
$\tilde{w}^k\notin\partial T(V\cdot)(0)$, i.e., there exists
\begin{equation}\label{wk-cond0}
w^k\in\partial\|Q^{\frac12}V\cdot\|(0)~\text{s.t.}~E(l^k)w^k\notin\partial T(V\cdot)(0).
\end{equation}
Based on \cite[\S 3.5]{Pham1998}, by selecting such a $w^k$,
we can guarantee that the computed $t^k$ fulfills the condition
$(\tilde{g}-\tilde{h})(t^k)<(\tilde{g}-\tilde{h})(0)=0$, i.e., $E(t^k)<E(l^k)$,
where $\tilde{g}$ and $\tilde{h}$ are defined as \eqref{param-rho}.
Condition \eqref{wk-cond0} refers to
$$
w^k\in V^\top Q^{\frac12}\partial\|\cdot\|(0)~~\text{s.t.}~~E(l^k)w^k\notin V^\top
\partial T(0).
$$
Note that selecting $w^k$ has been translated into the act of choosing $c$
from $\mathcal{B}(0,\sqrt{d_{min}})$, hence we need to
select $c\in\br^m$ to satisfy
\begin{equation}\label{c-condition}
c\in\mathcal{B}(0,\sqrt{d_{min}})~\text{s.t.}~ E(l^k)c\notin V^\top\partial T(0).
\end{equation}

For any given pair $(i,j)$, we can assume without loss of generality that
$w_{ij}\leq w_{ji}$,
since the case $w_{ij}\geq w_{ji}$ can be discussed similarly. Then it holds
\begin{equation*}
\begin{split}
w_{ij}[x_i-x_j]_++w_{ji}[x_j-x_i]_+&=w_{ij}[x_i-x_j]_++w_{ij}[x_j-x_i]_+
+(w_{ji}-w_{ij})[x_j-x_i]_+\\
&=w_{ij}([x_i-x_j]_++[x_j-x_i]_+)+(w_{ji}-w_{ij})[x_j-x_i]_+\\
&=w_{ij}|x_i-x_j|+(w_{ji}-w_{ij})[x_j-x_i]_+.
\end{split}
\end{equation*}
Denote $I_{min}=\{(i,j)\mid w_{ij}<w_{ji},1\leq i,j\leq n\}$,
$I_w=\{(i,j)\mid w_{ij}=w_{ji},1\leq i<j\leq n\}$. Then we can express $T(x)$ as
$$T(x)=\sum_{i,j=1}^nw_{ij}[x_i-x_j]_+=\sum_{(i,j)\in I_{min}\cup I_w}
w_{ij}|x_i-x_j|+(w_{ji}-w_{ij})[x_j-x_i]_+.
$$
Hence, by applying the conclusions drawn from \cite[Example 3.44, Theorem 3.50]
{Beckoptim2017}, we have
\begin{equation}\label{partial T0}
\partial T(0)=\sum_{(i,j)\in I_{min}\cup I_w}w_{ij}a_{ij}m_{ij}+(w_{ji}-w_{ij})
b_{ji}m_{ji},
\end{equation}
where $|a_{ij}|\leq1,~0\leq b_{ji}\leq1$, $m_{ij}\in\br^n$ satisfies
$m_{ij}(s)=0, \forall s\notin\{i,j\}$, and
$m_{ij}(i)=1,m_{ij}(j)=-1$.

Therefore, for any $\tilde{v}\in\partial T(0)$ and $s=1,2,\cdots,m$,
\begin{equation}\label{V-subgT}
(V^\top \tilde{v})_s=\sum_{(i,j)\in I_{min}\cup I_w}[w_{ij}(v_s(i)-v_s(j))a_{ij}+
(w_{ji}-w_{ij})(v_s(j)-v_s(i))b_{ji}],
\end{equation}
where $v_s$ is the $s$-th column of $V$.
From \eqref{V-subgT} we can deduce that
$$
-T_1(v_s)\leq(V^\top \tilde{v})_s\leq T(v_s),~\forall s,
$$
where $T_1(x)=\sum_{(i,j)\in I_{min}\cup I_w}
w_{ij}|x_j-x_i|+(w_{ji}-w_{ij})[x_i-x_j]_+=
\sum_{i,j=1}^nw_{ij}[x_j-x_i]_+$.
Consequently, $E(l^k)c\in V^\top\partial T(0)$ implies
\begin{equation}\label{bound-T-subdiff}
-T_1(v_s)\leq E(l^k)c_s\leq T(v_s),
\forall 1\leq s\leq m.
\end{equation}
 Denote $i_m=\argmin_{i=1,\cdots,m}T(v_i)$,
$\bar{i}_m=\argmin_{i=1,\cdots,m}T_1(v_i)$.
Combining \eqref{bound-T-subdiff} and \eqref{c-condition},
we conclude that when
\begin{equation}\label{choose w nece cond}
\sqrt{d_{min}}E(l^k)>\min\{T(v_{i_m}),T_1(v_{\bar{i}_m})\},
\end{equation}
then there exists $c\in\mathcal{B}(0,\sqrt{d_{min}})$ such that
$E(l^k)c$ does not meet the requirement of \eqref{bound-T-subdiff}. In this
case, we can select $c$ in the following manner.
\begin{itemize}
  \item $T(v_{i_m})\leq T_1(v_{\bar{i}_m})$: Taking $c=\sqrt{d_{min}}e_{i_m}$,
  where $e_i$ is the $m$-length column vector whose $i$th component is one while all the others are zeros.
  \item $T(v_{i_m})>T_1(v_{\bar{i}_m})$: Taking $c=-\sqrt{d_{min}}e_{\bar{i}_m}$.
\end{itemize}
If $l^k$ fails to meet \eqref{choose w nece cond}, i.e.,
$\sqrt{d_{min}}E(l^k)\leq \min\{T(v_{i_m}),T_1(v_{\bar{i}_m})\}$,
then for any $c\in\mathcal{B}(0,\sqrt{d_{min}})$, $E(l^k)c$ fulfills the necessary
condition \eqref{bound-T-subdiff}. However, there may still exist a number $c$
that belongs to $\mathcal{B}(0,\sqrt{d_{min}})$, while $E(l^k)c$ does not belong to
$V^\top\partial T(0)$. In such cases, we lack an effective method
to identify such $c$, so we have to randomly select one. The strategy
involves selecting an element from the set $C:=\pm\sqrt{d_{min}}
\{e_i\}_{i=1}^m$ during each iteration.
The random selection method, together
with the selection method under scenario \eqref{choose w nece cond}, is
summarized in Algorithm \ref{alg:choose w}. We remark that if the graph is
undirected, the set $\partial T(0)$ is symmetric, i.e., $-x\in\partial T(0)$
can be inferred from $x\in\partial T(0)$. Therefore, whether $E(l^k)w^k$ belongs
to $V^\top\partial T(0)$ is independent of its sign, whereas in the case of
directed graph, $\partial T(0)$ is generally asymmetric and
the choice of different signs for $w^k$ may influence whether $E(l^k)w^k$
falls within $V^\top\partial T(0)$.
In concluding this section, we
would like to point out that when $l^k$
meets \eqref{choose w nece cond}, selecting $w^k$
in accordance with Algorithm \ref{alg:choose w} can guarantee the improvement
of the candidate point $l^k$. According to our experiment, it is primarily the
first selection method that proves effectiveness (see Section \ref{exp}
for details).

\begin{algorithm}[H]
\caption{The method of Selecting $w^k$ in PS-DCA for Solving Problem \eqref{GFM-fracprogram}}
\label{alg:choose w}
\begin{algorithmic}[1]
\State Input $\sqrt{d_{min}}$, $i_m$ and $T(i_m)$, $\bar{i}_m$ and $T(\bar{i}_m)$.
 \If{$\sqrt{d_{min}}E(l^k)>\min\{T(v_{i_m}),T_1(v_{\bar{i}_m})\}$}
 \State if $T(v_{i_m})\leq T_1(v_{\bar{i}_m})$, set $w^k=\sqrt{d_{min}}e_{i_m}$. Otherwise set $w^k=-\sqrt{d_{min}}e_{\bar{i}_m}$.
 \Else
 \State randomly selecting $i\in\{1,\cdots,m\}$, set
 $w^k=\pm\sqrt{d_{min}}e_i$, in which the positive and negative signs are
 randomly allocated.
 \EndIf
\end{algorithmic}
\end{algorithm}

\section{Numerical Experiments}
\label{exp}
In this section, we conduct some numerical experiments to test the efficiency
of the proposed PS-DCA. We concentrate on calculating the generalized graph
Fourier modes
given by Definition \ref{de-GFM-gdv} with $Q=I$ and $Q=D$,
where $D$ denotes the diagonal degree matrix of
the graph.
We compute these generalized graph Fourier modes on three different types of graphs:
community graph, random geometric graph (RGG) and
directed random geometric graph (DRGG), where the community graph and RGG are undirected,
while the DRGG is a directed version of the RGG. The community graph is produced using GSPBOX\cite{Perraudin2014GSPBOXAT}.
The RGG is generated by $n$ random points $p_i$ chosen uniformly in the unit
square $[0,1]^2$, and the weights are defined by
$w_{ij}:=\exp(-\|p_i-p_j\|^2/0.5)$.
To prevent the RGG from becoming overly dense, we reset the weights to $0$ if they are less than $0.7$.
By making slight modifications to the generation method of RGG,
we obtain the DRGG. More specifically,
we firstly generate $n$ random points $p_i$ uniformly distributed within
the square $[0,1]^2$, and the weight $w_{ij}$ connecting $p_i$ and $p_j$ is set to
$1$ with probability $1-\exp(-\|p_i-p_j\|^2/0.5)$, or set to $0$ otherwise.
Both the RGG and DRGG are generated with the help of GSPBOX.

We use the notation $\mathcal{G}=(\mathcal{V},\mathcal{E})$ to represent
a weighted graph, where $\mathcal{V}=\{1,2,\cdots,n\}$ denotes its vertex set,
and $\mathcal{E}=\{w_{ij}\}_{i,j\in\mathcal{V}}$ denotes its edge set
satisfying $w_{ij}>0$ if there is a edge from vertex $i$ to vertex $j$
and $w_{ij}=0$ otherwise.
Let $|\mathcal{E}|$ be the cardinality of $\mathcal{E}$.
The number of vertices and edges of the graphs in our experiment
is specified in Table \ref{size of graphs}. In DRGG,
the number of directed edges makes up roughly 49\% of the total edge count.
All experiments were carried out in Matlab on a laptop
with an Intel(R) Core(TM) i7-1165G7 CPU (2.80 GHz) and 16 GB of RAM.
\begin{table}[h]
\centering
\caption{The number of vertices and edges for the tested graphs}
\label{size of graphs}
\resizebox{0.5\textwidth}{!}{
\begin{tabular}{cccc}
\toprule
\multicolumn{1}{c}{}&\textbf{community}&\textbf{RGG}&\textbf{DRGG}\\
\midrule[.5pt]
\multicolumn{1}{c}{$(n,|\mathcal{E}|)$}&$(20,97)$&$(20,75)$&$(20,143)$\\
\bottomrule
\end{tabular}}
\end{table}

Since the calculation of each (T,Q)-GFM can be equivalently reduced to
solving problem \eqref{GFM-fracprogram}, in the cases of $Q=I$ and $Q=D$, we
apply PS-DCA, PGSA, and PGSA\_BE to solve problem \eqref{GFM-fracprogram}
for given $k$ and $U_{k-1}$, with all algorithms starting from the same
initial point.
When $Q=I$,
we employ the graph Laplacian matrix to generate
$U_{k-1}$. Specifically, denoted by $U$ the Laplacian basis of
the graph, which can be computed using the function
``G=gsp\_compute\_fourier\_basis(G)'' in GSPBOX.
Let the first $k-1$ column vectors of $U$ be $U_{k-1}$.
Note that when the graph is directed, GSPBOX also offers an available type of
Laplacian matrix for utilization. For a detailed description, please refer to
\cite[Chapter 9]{Perraudin2014GSPBOXAT}.
When $Q=D$, we begin by generating
a matrix $U$ such that $U^\top QU=I$, where the first column of $U$ is identical
to $u_1$ as specified in Definition \ref{de-GFM-gdv}, and then designate the
first $k-1$ columns of $U$ as $U_{k-1}$.
$U$ is generated in the following manner: First, set $v_1=Q^{\frac12}u_1$. Then,
utilize Matlab function ``null($v_1^\top$)'' to produce a set of orthonormal
vectors orthogonal to $v_1$, denoted as $\tilde{V}$. Finally,
let $U = Q^{-\frac12}[v_1,\tilde{V}]$.  For both of the two cases,
the matrix $V\in\br^{n\times m}$ in \eqref{GFM-fracprogram} is
produced by Matlab function ``null($U_{k-1}^\top Q$)''.

\begin{table}[t]
\centering
\caption{Computational results of the $k$-th (T,I)-GFM averaged over the first 30\% initial points.}
\label{T-I-GFMfirst}
\resizebox{\textwidth}{!}{
\begin{tabular}{cccccccccc}
\toprule
\multicolumn{1}{c}{\multirow{2}{*}{}}&\multicolumn{4}{c}{Rand.Geometric}&
& \multicolumn{4}{c}{Directed Rand.Geometric}\\
\cmidrule{2-5}\cmidrule{7-10}
\multicolumn{1}{c}{}&\textbf{PS-DCA}&\textbf{PGSA}&\textbf{PGSA\_BE}&\textbf{ManPG-Ada}
&&\textbf{PS-DCA}&\textbf{PGSA}&\textbf{PGSA\_BE}&\textbf{ManPG-Ada}\\
\midrule[.5pt]
\multicolumn{1}{c}{$\bf{k=2}$}\\
\multicolumn{1}{c}{Obj}&1.8581&9.4254&9.4254&9.4254&&5.1299&14.1585&14.1585
&14.1585\\
\multicolumn{1}{c}{Iter}&2&2.2667&2.1333&2&&
2&3.1333&2.8667&2\\
\multicolumn{1}{c}{Time}&1.6857&1.015&0.9485&0.6898&&
1.7542&1.5133&1.3629&0.7879\\
\multicolumn{1}{c}{$\bf{k=3}$}\\
\multicolumn{1}{c}{Obj}&3.1521&9.7009&9.7009&9.7009&&8.1352&
15.1768&15.1768&15.1914\\
\multicolumn{1}{c}{Iter}&2.2&2.3333&2.3333&2&&3&3
&3&2.4\\
\multicolumn{1}{c}{CPU time}&1.8978&1.0881&1.0832&0.6848&&
2.6425&1.5009&1.5072&0.9639\\
\multicolumn{1}{c}{$\bf{k=4}$}\\
\multicolumn{1}{c}{Obj}&3.2709&9.967&9.967&9.967
&&10.0946&15.0618&15.0618&15.0846\\
\multicolumn{1}{c}{Iter}&3&4&4&3&&
2.2&3&3&2.5333\\
\multicolumn{1}{c}{CPU time}&2.6015&1.8897&1.9071&1.083&&
1.9559&1.4376&1.4694&1.0373\\
\multicolumn{1}{c}{$\bf{k=5}$}\\
\multicolumn{1}{c}{Obj}&6.3784&10.3004&10.3004&10.3004&&10.7268
&15.4761&15.4761&15.6138\\
\multicolumn{1}{c}{Iter}&3&3.4667&3.4667&3&&2&3.2667&
3.2667&2.6667\\
\multicolumn{1}{c}{CPU time}&2.7098&1.7419&1.8065&1.1318
&&1.7379&1.5509&1.5726&1.0545\\
\bottomrule
\end{tabular}}
\end{table}

Generally speaking, the performance of algorithms designed to solve
nonconvex problems is significantly influenced by the choice of the initial point.
To examine the impact of initial points on the algorithm's computation of each
generalized graph Fourier mode,
we ran PS-DCA, PGSA, and PGSA\_BE with 50 initial points for each $k$.
The first 30\% of the initial points are linear combinations of two column
vectors that have the two largest objective function values from the remaining
$n-(k-1)$ columns of $U$, with the combination coefficients being randomly and
uniformly selected from $(0,1)$. While the remaining 70\%
of the initial points are acquired by projecting randomly generated points onto
$\mathcal{S}$. The step size $\lambda_k$ employed in the proximal subgradient
iteration \eqref{pgsa} proposed by \cite{PGSA2022} and \cite{PGSABE2022} is
bounded by the Lipschitz constant associated with $h$ in \eqref{structured FP},
for problem \eqref{GFM-fracprogram}, however, this constant can be chosen as any
positive number.
Inspired by the stepsize setting in \cite[\S 5]{smoothalg2024}
for evaluating the proximal operator of the graph directed total variation
at a point,
we tune this parameter and finally set $\lambda_k=80/\|CV\|$ for PGSA and $\alpha=40/\|CV\|$ for PGSA\_BE,
where $C\in\br^{|\mathcal{E}|\times n}$ encodes the connectivity of the graph,
with the $k$-th row of $C$ corresponding to the $k$-th directed link $w_{i_k,j_k}$
of the graph such that $C_{ki_k}=1$ and $C_{kj_k} = -1$, while all other elements
of $C$ are zero. In addition, for the remaining parameters of PGSA\_BE, we adopt those
recommended in \cite{PGSABE2022}: $\epsilon=10^{-4}$,
$\beta_k=(\theta_{k-1}-1)/\theta_k$, where $\theta_{-1}=\theta_0=1$,
$\theta_{k+1}=(1+\sqrt{1+4\theta^2_k})/2$, and we reset $\theta_k=1$ every
$10$ iterations. For the interpretation of these parameters,
readers may refer to page 105 of \cite{PGSABE2022}. For PS-DCA,
we set the acceptance parameter $\epsilon'=10^{-6}$.
The choice of regularization parameter $\rho_k$
in \eqref{param-rho} depends on the d.c. decomposition employed by the d.c.
function. Unfortunately, there lacks systematic research on how to tune this
parameter for a given d.c. decomposition.
We try our best to tune the proximal stepsize and the regularization parameter in PS-DCA,
and the ones we adopt in our experiment are $\lambda_k=100/\|CV\|$ and
$\rho_k=E(l^k)$. These three algorithms are terminated when
$|E(x^{k+1})-E(x^k)|<10^{-6}$ occurs or the iteration count reaches $20$.

\begin{table}[t]
\centering
\caption{Computational results of the $k$-th (T,I)-GFM averaged over all of the initial points.}
\label{T-I-GFMall}
\resizebox{\textwidth}{!}{
\begin{tabular}{cccccccccc}
\toprule
\multicolumn{1}{c}{\multirow{2}{*}{}}&\multicolumn{4}{c}{Rand.Geometric}&
& \multicolumn{4}{c}{Directed Rand.Geometric}\\
\cmidrule{2-5}\cmidrule{7-10}
\multicolumn{1}{c}{}&\textbf{PS-DCA}&\textbf{PGSA}&\textbf{PGSA\_BE}&\textbf{ManPG-Ada}
&&\textbf{PS-DCA}&\textbf{PGSA}&\textbf{PGSA\_BE}&\textbf{ManPG-Ada}\\
\midrule[.5pt]
\multicolumn{1}{c}{$\bf{k=2}$}\\
\multicolumn{1}{c}{Obj}&1.8567&4.4579&4.4006&4.4612&&5.1299&9.9015&9.9015
&9.9758\\
\multicolumn{1}{c}{Iter}&2.04&6.52&6.28&2.5&&
2&3.96&3.78&2.18\\
\multicolumn{1}{c}{Time}&1.8516&2.8994&2.8141&0.9313&&
1.763&1.904&1.7957&0.8626\\
\multicolumn{1}{c}{$\bf{k=3}$}\\
\multicolumn{1}{c}{Obj}&3.1521&5.9743&5.9743&5.9461&&8.1569&
11.9653&11.9653&12.3278\\
\multicolumn{1}{c}{Iter}&2.12&4.64&4.54&2.48&&2.96&4.66
&4.58&2.66\\
\multicolumn{1}{c}{CPU time}&1.9101&2.193&2.1494&0.9235&&
2.7728&2.4383&2.4202&1.1295\\
\multicolumn{1}{c}{$\bf{k=4}$}\\
\multicolumn{1}{c}{Obj}&3.2709&6.7208&6.7208&6.7512
&&10.1846&13.5711&13.5711&13.8282\\
\multicolumn{1}{c}{Iter}&2.88&4.88&4.92&2.98&&
2.5&4.08&4.04&2.84\\
\multicolumn{1}{c}{CPU time}&2.5507&2.3446&2.76&1.1337&&
2.2033&1.9541&1.9975&1.1518\\
\multicolumn{1}{c}{$\bf{k=5}$}\\
\multicolumn{1}{c}{Obj}&6.4197&8.0319&8.031&7.9608&&10.7268
&14.2687&14.2687&14.3697\\
\multicolumn{1}{c}{Iter}&3.16&4.96&5.14&3.18&&2.1&4.08&
4.06&3.06\\
\multicolumn{1}{c}{CPU time}&3.1206&2.6922&2.8685&1.3347
&&1.9214&2.0923&2.1522&1.2942\\
\bottomrule
\end{tabular}}
\end{table}

For the purpose of comparing it with algorithms for solving fractional
minimization problems, we also utilize ManPG-Ada (manifold proximal gradient
method incorporating adaptive ticks)\cite{manpg2020} to solve the original
problem \eqref{de-T-Q-GFM} with the same $U_{k-1}$ as the aforementioned three
algorithms. The ManPG is designed to tackle nonsmooth optimization problems over
the Stiefel manifold. We reformulate \eqref{de-T-Q-GFM} as the following
spherical constrained optimization problem:
\begin{equation}\label{equiva T-Q}
\begin{split}
\min_{x\in\br^m}&~T(Q^{-\frac12}\tilde{V}x)\\
        s.~t. &~\|x\|=1,
\end{split}
\end{equation}
where $\tilde{V}=(\tilde{v}_1,\cdots,\tilde{v}_m)$ with $\{\tilde{v}_i\}^m_{i=1}$ is an orthonormal basis for the space $\{x\in\br^n\mid U^\top_{k-1}Q^{\frac12}x=0\}$.
When $Q=I$, both $\tilde{V}$ in \eqref{equiva T-Q} and $V$ in \eqref{GFM-fracprogram} are formed by the orthonormal basis of space $\{x\in\br^n\mid U^\top_{k-1}x=0\}$. Hence, we set $\tilde{V}=V$
and initialize ManPG-Ada at the same point as PS-DCA, PGSA, and PGSA\_BE.
In the case of $Q=D$, the $\tilde{V}$ in \eqref{equiva T-Q} is produced in Matlab by ``null($U_{k-1}^\top D^{\frac12}$)'' and the
ManPG-Ada is initialized at the point that shares the same initial objective function value as that in the fractional programming algorithms being compared.
Note that when ManPG-Ada is applied to solve \eqref{equiva T-Q},
the primary computational cost in each iteration lies on the
nonsmooth convex subproblem for finding the descent direction:
\begin{equation}\label{subp-ManPG}
v^i:=\argmin_{v\in\br^m} \tilde{T}(x^{i}+v)+\frac1{2t}\|v\|^2,
~~\text{s.t.}~v^\top x^i=0
\end{equation}
where
$\tilde{T}(x)=T(Q^{-\frac12}\tilde{V}x)$, $t>0$ is the stepsize.
Since the proximal operator of $\tilde{T}$ has no closed-form solution,
the semi-smooth Newton method,
which is recommended in \cite{manpg2020}, cannot be utilized for efficiently solving
subproblem \eqref{subp-ManPG}. As an alternative, we use the convex optimization
toolbox CVX\cite{cvx} with the default solver SDPT3 to solve it. Analogous to the
implementation of computing the Fourier basis for a directed graph utilizing
ManPG-Ada\cite{smoothalg2024}, we set $t=100/\|BQ^{-\frac12}\tilde{V}\|$,
$\tau=1.01$, $\gamma=0.5$. For the meaning of these parameters, we refer the
readers to page 226-227 of \cite{manpg2020}.
Moreover, ManPG-Ada is terminated when $\|v^i/t\|^2<10^{-8}n^2$
or its iteration number hits $20$.

\begin{table}[t]
\centering
\caption{Computational results of the $k$-th (T,D)-GFM averaged over the first
30\% initial points.}
\label{T-D-GFMfirst}
\resizebox{\textwidth}{!}{
\begin{tabular}{cccccccccc}
\toprule
\multicolumn{1}{c}{\multirow{2}{*}{}}&\multicolumn{4}{c}{Community}&
& \multicolumn{4}{c}{Directed Rand.Geometric}\\
\cmidrule{2-5}\cmidrule{7-10}
\multicolumn{1}{c}{}&\textbf{PS-DCA}&\textbf{PGSA}&\textbf{PGSA\_BE}&\textbf{ManPG-Ada}
&&\textbf{PS-DCA}&\textbf{PGSA}&\textbf{PGSA\_BE}&\textbf{ManPG-Ada}\\
\midrule[.5pt]
\multicolumn{1}{c}{$\bf{k=2}$}\\
\multicolumn{1}{c}{Obj}&0.7554&2.562&2.562&2.0435&&1.6191&4.0912&4.0912
&4.0912\\
\multicolumn{1}{c}{Iter}&4&3.2667&3.2667&2.2&&
2&2.1333&2.0667&2\\
\multicolumn{1}{c}{Time}&3.5085&1.6034&1.6351&0.8442&&
1.7902&1.0232&0.9853&0.7618\\
\multicolumn{1}{c}{$\bf{k=3}$}\\
\multicolumn{1}{c}{Obj}&2.257&3.7306&3.7306&3.744&&1.6958&
4.1169&4.1169&4.1434\\
\multicolumn{1}{c}{Iter}&2&3.4&3.3333&2.0667&&2&2.4667
&2.3333&2\\
\multicolumn{1}{c}{CPU time}&1.756&1.6532&1.598&0.7491&&
1.8745&1.1926&1.1969&0.8497\\
\multicolumn{1}{c}{$\bf{k=4}$}\\
\multicolumn{1}{c}{Obj}&2.257&4.0066&4.0066&3.7878
&&1.8775&4.0833&4.0833&4.0833\\
\multicolumn{1}{c}{Iter}&2&2.6&2.5333&2.2&&
2&2.2&2.2&2\\
\multicolumn{1}{c}{CPU time}&1.978&1.4197&1.475&0.9009&&
1.9154&1.1567&1.1028&0.7955\\
\multicolumn{1}{c}{$\bf{k=5}$}\\
\multicolumn{1}{c}{Obj}&2.257&3.9919&3.9919&3.7731&&2.0188
&4.2121&4.2121&4.2121\\
\multicolumn{1}{c}{Iter}&2&2.7333&2.7333&2.2&&2&2.7333&
2.7333&2.3333\\
\multicolumn{1}{c}{CPU time}&1.892&1.4602&1.4439&0.8623
&&1.8465&1.4118&1.3548&0.9065\\
\bottomrule
\end{tabular}}
\end{table}

\begin{table}[t]
\centering
\caption{Computational results of the $k$-th (T,D)-GFM averaged over all of the initial points.}
\label{T-D-GFMall}
\resizebox{\textwidth}{!}{
\begin{tabular}{cccccccccc}
\toprule
\multicolumn{1}{c}{\multirow{2}{*}{}}&\multicolumn{4}{c}{Community}&
& \multicolumn{4}{c}{Directed Rand.Geometric}\\
\cmidrule{2-5}\cmidrule{7-10}
\multicolumn{1}{c}{}&\textbf{PS-DCA}&\textbf{PGSA}&\textbf{PGSA\_BE}&\textbf{ManPG-Ada}
&&\textbf{PS-DCA}&\textbf{PGSA}&\textbf{PGSA\_BE}&\textbf{ManPG-Ada}\\
\midrule[.5pt]
\multicolumn{1}{c}{$\bf{k=2}$}\\
\multicolumn{1}{c}{Obj}&0.7554&1.7768&1.6259&1.6368&&1.6295&2.9264&2.9158
&2.9371\\
\multicolumn{1}{c}{Iter}&3.92&3.92&4.24&2.26&&
2.04&3.12&3.24&2.1\\
\multicolumn{1}{c}{Time}&3.732&2.0958&2.3054&0.8911&&
1.9766&1.622&1.6585&0.9084\\
\multicolumn{1}{c}{$\bf{k=3}$}\\
\multicolumn{1}{c}{Obj}&2.1908&2.8915&2.8915&2.8327&&1.6982&
3.3548&3.3753&3.3648\\
\multicolumn{1}{c}{Iter}&2.04&4.3&4.16&2.26&&2&3.24
&3.24&2.12\\
\multicolumn{1}{c}{CPU time}&1.9564&2.3113&2.1682&0.8696&&
1.8722&1.6444&1.639&0.93\\
\multicolumn{1}{c}{$\bf{k=4}$}\\
\multicolumn{1}{c}{Obj}&2.257&3.7586&3.7586&3.561
&&1.8862&3.4274&3.4274&3.4119\\
\multicolumn{1}{c}{Iter}&2&3.42&3.36&2.28&&
2.02&3.62&3.66&2.32\\
\multicolumn{1}{c}{CPU time}&1.9803&1.9289&1.8968&0.9184&&
2.04&2.0432&2.051&1.0813\\
\multicolumn{1}{c}{$\bf{k=5}$}\\
\multicolumn{1}{c}{Obj}&2.257&3.7843&3.7843&3.6953&&2.0223
&3.4448&3.4448&3.4345\\
\multicolumn{1}{c}{Iter}&2.02&3.16&3.18&2.2&&2.02&3.72&
3.82&2.6\\
\multicolumn{1}{c}{CPU time}&1.9651&1.7521&1.7434&0.8765
&&1.857&1.9067&1.9846&1.082\\
\bottomrule
\end{tabular}}
\end{table}

We compute the first five modes for (T,I)-GFMs and (T,D)-GFMs for two different types
of graphs, with the results displayed in Tables \ref{T-I-GFMfirst} to \ref{T-D-GFMall}.
Since the first modes of both (T,I)-GFMs and (T,D)-GFMs correspond to
constant graph signals, we begin with the calculation of the second
generalized graph Fourier mode. Table \ref{T-I-GFMfirst} (Table \ref{T-D-GFMfirst})
reports the mean objective value,
iteration number and CPU time (in seconds) obtained from the computation of (T,I)-GFMs ((T,D)-GFMs),
averaged over the first 30\% initial points,
which typically possess comparatively larger objective function values,
while Table \ref{T-I-GFMall} (Table \ref{T-D-GFMall})
reports the computational results of (T,I)-GFMs ((T,D)-GFMs)
averaged over all of the initial points. From Table \ref{T-I-GFMall} and Table \ref{T-D-GFMall}
we observe that the objective values of the solutions provided by PGSA, PGSA\_BE
and ManPG-Ada are very close but notably higher than that obtained by PS-DCA.
The reason lies in the fact that PS-DCA incorporates an additional step of DCA
starting from $0$ during the iteration, preventing it from getting trapped in
poor-quality local minimizers, unlike the other three algorithms.
Furthermore, compared to PGSA and PGSA\_BE, PS-DCA requires fewer iterations,
comparable to ManPG-Ada, but the CPU time required is comparable to or
slightly higher than that of PGSA and PGSA\_BE.
This is because during the iteration, if DCA identifies a point superior to the
candidate point, PS-DCA will commence reducing the objective function value
from this new point, potentially increasing the number of iterations required by
the algorithm. Additionally,
in some instances, DCA's
failure to enhance the candidate point can further increase
the algorithm's running time. We also notice that the CPU time required by these
three fractional programming algorithms exceeds that of ManPG-Ada,
which can be viewed as performing the steepest descent method on
$\tilde{T}$ within the Stiefel manifold.

From Table \ref{T-I-GFMfirst} and Table \ref{T-D-GFMfirst} we see that
in the computation of all graph Fourier modes,
when the initial point is inadequately selected, PGSA, PGSA\_BE, and ManPG-Ada are
unable to find a satisfactory solution with respect to the
graph directed variation of the signal. By comparing Table \ref{T-I-GFMfirst} with
Table \ref{T-I-GFMall}, we observe that for each $k$,
the influence of initial points on PS-DCA,
in terms of objective function value, is significantly lower than that of
the other algorithms, and the impact of initial points on PS-DCA is also
less pronounced than that on PGSA and PGSA\_BE in terms of iteration number
and CPU time.
A similar conclusion can be drawn for (T,D)-GFMs
by comparing Table \ref{T-D-GFMfirst} with Table \ref{T-D-GFMall}.
This demonstrates that introducing the one-step DCA can substantially
diminish the algorithm's reliance on initial points,
thereby improving the stability of the algorithm's performance.

\section{Conclusion}
\label{conclud}
In this paper, we study a class of fractional minimization
problems characterized by a single ratio, in which both the numerator and
denominator of the objective are convex functions that exhibit positive
homogeneity. By applying convex analysis theory,
we derive the local optimality conditions pertaining to this problem,
as well as an equivalent condition for obtaining a global optimal solution.
Then we propose a proximal algorithm incorporating d.c. optimization, namely,
PS-DCA to find the critical points of the problem. When the DCA step within PS-DCA
is removed, it reduces to the PSA.
We derive subsequential convergence property of PSA or PS-DCA
under mild assumptions and establish closedness property with respect to local
minimizers of $E$ for PSA. Note that, generally, PS-DCA is not closed with
respect to local minimizers of $E$, particularly when the candidate point is
located near a local minimizer with a relatively high objective function value.
Thanks to our subgradient selection strategy in the DCA step, subsequent iterates
will escape the neighborhood of this local minimizer and move to a region with
considerably lower objective values. Moreover,
for a typical class of generalized graph Fourier mode problems, we establish
global convergence of PSA or PS-DCA under
reasonable assumption regarding the zeros of $T(x)$.
Numerical experiments of computing the generalized graph Fourier modes
demonstrate that,
in comparison to existing methods, the proposed PS-DCA yields
notably smaller graph directed variation and exhibits robustness in its
performance, particularly when confronted with improperly chosen initial points,
it remains capable of computing valid graph Fourier modes.

It is worth noting that, according to Theorem \ref{simp dca converg thm} and
Theorem \ref{equ cond P3}, the quality of the
solution provided by PS-DCA depends on whether DCA can find a point
superior to the candidate, which ultimately boils down to the selection of the
subgradient $w^k$. In subsection \ref{choos-wk}, we discuss the selection
method of $w^k$ for a class of generalized graph Fourier mode problems, and
our experiments further validate the efficiency of this method.
For future improvement directions,
one may consider the construction of $V$, which incorporates columns that
result in a relatively smaller $T(x)$ or $T_1(x)$. Given that many different
forms of graph variations have been proposed so far,
developing the corresponding selection method for $w^k$ remains an
interesting and meaningful problem.

\section*{Acknowledgement}
This research was partially supported by
National Natural Science Foundation of China (No: 12171488),
and the Talent Special Project of Guangdong Polytechnic
Normal University (Grants No. 2021SDKYA083)

\begin{appendices}

\section{Proof of the equivalence between problem \eqref{eq:P1} and
fractional programming \eqref{eq:P2}}
\label{appendix}
It is easy to verify that positively homogeneous convex functions must
satisfy triangle inequality:
\begin{equation}\label{eq:triangle-inequality}
T(x+y)\le T(x)+T(y),~~~B(x+y)\le B(x)+B(y),
~~~~\forall x, y\in X.
\end{equation}
In fact, by utilizing the convexity of $T$, we have
$$T(x+y)=2T\Big(\frac{1}{2}x+\frac{1}{2}y\Big)\le T(x)+T(y),$$
and similarly the second inequality of \eqref{eq:triangle-inequality}. Hence,
$B(x)$ is a seminorm, which is not true for $T(x)$ because it only satisfies
positive homogeneity\cite{Yosida-k1995}.

For any $x, y\in X_0$ and $a,b\in\mathbb{R}$, by using triangle inequality,
we obtain
$$0\le B(ax+by)\le B(ax)+B(by)=|a|B(x)+|b|B(y)=0,$$
which implies that $X_0$ is a linear subspace of $X$. Obviously, we only need to
solve problem \eqref{eq:P1} in the orthogonal complementary space $X_0^\perp$.
Indeed, the restraint of $B$ on subspace $X_0^\perp$ satisfies $B(x)=0$ if and
only if $x=0$ and thus $B$ is a norm on $X_0^\perp$.
For any $x\in X$, we have $x=y+z$, where $y\in X_0,~z\in X_0^\perp$, then
$$T(z)=T(x-y)\le T(x)+T(-y)=T(x),~~~~T(x)\le T(y)+T(z)=T(z).$$
Therefore $T(x)=T(z)$. Likewise we can prove $B(x)=B(z)$.
As a result, $x$ is the solution to problem \eqref{eq:P1} if and only if
$x=y+z$, where $y\in X_0$, $z$ is the solution to the following problem
$$\begin{cases}
\min_{z\in X_0^\perp}T(z),\\
\text{s.t.}~~ B(z)=1.
\end{cases}$$
\par
We take an orthonormal basis $v_1,\cdots,v_m$ of $X_0^\perp$.
Let $V:=(v_1,\cdots,v_m)$, then for any $z\in X_0^\perp$, we have $z=Vx$ for
some  $x\in\br^m$, thus problem \eqref{eq:P1} boils down to
\begin{equation}\label{eq:P1-1}
\begin{cases}
\min_{x\in\br^m}T(Vx) \\
\text{s.t.}~~B(Vx)=1,
\end{cases}
\end{equation}
where $B(Vx)>0,~\forall x\in\br^m\setminus\{0\}$.
\par
For the sake of simplicity, we still use $T(x)$ and $B(x)$ to denote $T(Vx)$
and $B(Vx)$ respectively, thus problem \eqref{eq:P1-1} can be written as
\begin{equation*}\label{eq:tilP1}
\tilde{P}:~~\begin{cases}
\min_{x\in\br^m}T(x) \\
\text{s.t.}~~ B(x)=1,
\end{cases}
\end{equation*}
where $T$ and $B$ are convex functions on $\br^m$ satisfying positive
homogeneity and absolute homogeneity, respectively, and
\begin{equation}\label{eq:normalization}
B(x)>0,~~~~\forall x\in \br^m\setminus\{0\}.
\end{equation}
\par
The following lemma indicates that when condition \eqref{eq:normalization}
is satisfied, problem \eqref{eq:P1} is equivalent to fractional
programming \eqref{eq:P2}.

\begin{prop}\label{specP1-P2}
Let $x^*\in\br^m$. If $x^*$ is a (local) minimizer of Problem $\tilde{P}$,
then $x^*$ is also a (local) minimizer of problem \eqref{eq:P2}; On the contrary,
if $x^*$ is a (local) minimizer of problem \eqref{eq:P2}, then
$\bar{x}:=x^*/B(x^*)$ is a (local) minimizer of Problem $\tilde{P}$.
\end{prop}
\prf
We only provide proof for the local minimum case, the proof of global minimum is similar and simpler.
\par
Let $x^*$ be a local minimizer of Problem $\tilde{P}$, then $B(x^*)=1$,
and there exists $\delta>0$, such that
$$T(x^*)\le T(x),~~~~\forall x\in\{x\in\br^n|B(x)=1,~\|x-x^*\|<\delta\}.$$
\par
For any $x\in\br^n\setminus\{0\}$, $x\to x^*$ implies $B(x)\to B(x^*)=1$.
Therefore
$$\Big\|\frac{x}{B(x)}-x^*\Big\|
\le \|x-x^*\|+\Big|\frac{1}{B(x)}-1\Big|\|x\|\to 0.$$
Since $B(x/B(x))=1$, we have, when $\|x-x^*\|$ is sufficiently small,
$$\frac{T(x^*)}{B(x^*)}=T(x^*)
\le T\Big(\frac{x}{B(x)}\Big)=\frac{T(x)}{B(x)},$$
i.e., $x^*$ is a local minimizer of problem \eqref{eq:P2}.
\par
On the contrary, let $x^*\in\br^n$ be a local minimizer of problem \eqref{eq:P2},
then $x^*\neq 0$, and there exists $\delta>0 $, such that
$$\frac{T(x^*)}{B(x^*)}\le \frac{T(x)}{B(x)},~~~~
\forall x\in\br^n\setminus\{0\},~~~\|x-x^*\|<\delta.$$
Denote $\bar{x}:=x^*/B(x^*)$, then $B(\bar{x})=1$.
For $\forall x\in\br^n,~B(x)=1$, when $\|x-\bar{x}\|<\delta/B(x^*)$,
it holds $\|B(x^*)x-x^*\|<\delta$ and $B(x^*)x\neq 0$. Hence,
$$T(\bar{x})=\frac{T(x^*)}{B(x^*)}\le\frac{T(B(x^*)x)}{B(B(x^*)x)}=\frac{T(x)}{B(x)}
=T(x),$$
i.e., $\bar{x}$ is a local minimizer of Problem $\tilde{P}$.
\bbox

\section{Proof of Proposition \ref{loc trans prop}}
\label{append B}
\prf
It suffices to prove that there exists $\epsilon'>0$ such that $g(x)-h(x)\geq0,
\forall x\in\mathcal{B}(x^*,\epsilon')$. Firstly,
it is easy to verify that if there exists $\epsilon>0$ such that
$E(x^*)\le E(x),~\forall x\in\mathcal{B}_
{\scriptscriptstyle\mathcal{S}\cap\Omega}(x^*,\epsilon)$, if and only if
\begin{equation}\label{loc-E}
g(x)-h(x)\geq0,~~~~\forall x\in\mathcal{B}(x^*,\epsilon)\cap\mathcal{S}.
\end{equation}
Let $\epsilon'=\frac\epsilon{2+\epsilon}$, then for any
$x\in\mathcal{B}(x^*,\epsilon')$, we have
\begin{equation*}
\begin{split}
\|\frac{x}{\|x\|}-x^*\|&=\|\frac{x}{\|x\|}-\frac{x^*}{\|x^*\|}\|
\leq\|\frac{x}{\|x\|}-\frac{x^*}{\|x\|}\|+\|\frac{x^*}{\|x\|}-
\frac{x^*}{\|x^*\|}\|\\
&\leq\frac1{\|x\|}\epsilon'+|\frac1{\|x\|}-1|=\frac1{\|x\|}\frac\epsilon{2+\epsilon}
+|\frac1{\|x\|}-1|.
\end{split}
\end{equation*}
On the other hand, $x\in\mathcal{B}(x^*,\epsilon')$ implies that
$1-\epsilon'\leq\|x\|\leq1+\epsilon'$, i.e., $\frac2{2+\epsilon}\leq\|x\|\leq
\frac{2(1+\epsilon)}{2+\epsilon}$. Hence,
$$\|\frac{x}{\|x\|}-x^*\|\leq\frac{2+\epsilon}2\cdot\frac{\epsilon}{2+\epsilon}+
\frac{\epsilon}2=\epsilon.
$$
Consequently, $\frac x{\|x\|}\in\mathcal{B}(x^*,\epsilon)\cap\mathcal{S}$.
According to \eqref{loc-E}, it holds $g(\frac x{\|x\|})-h(\frac x{\|x\|})\geq0$,
i.e., $g(x)-h(x)\geq0$. The proof is complete.
\bbox

\end{appendices}

\bibliography{Myrefs}
\bibliographystyle{abbrv}

\end{document}